\documentclass[11pt]{amsart}

\newtheorem{theorem}{Theorem}[section]
\newtheorem{lemma}[theorem]{Lemma}
\newtheorem{proposition}[theorem]{Proposition}

\newtheorem{example}[theorem]{Example}
\newtheorem{definition}[theorem]{Definition}
\newtheorem{corollary}[theorem]{Corollary}

\newcommand{\N}{\mbox{$\mathbb{N}$}}

\newcommand{\Di}{\mbox{$\mathbb{D}$}}

\newcommand{\A}{\mbox{${\mathcal A}$}}
\newcommand{\B}{\mbox{${\mathcal B}$}}
\newcommand{\C}{\mbox{${\mathcal C}$}}
\newcommand{\D}{\mbox{${\mathcal D}$}}
\newcommand{\E}{\mbox{${\mathcal E}$}}
\newcommand{\F}{\mbox{${\mathcal F}$}}
\newcommand{\G}{\mbox{${\mathcal G}$}}
\newcommand{\Li}{\mbox{${\mathcal L}$}}
\newcommand{\U}{\mbox{${\mathcal U}$}}

\begin{document}
\title[Modules over operator algebras, and C$^*-$dilations]
{Modules over operator algebras, \\
 and the maximal C$^*-$dilation}
 
\vspace{30 mm}
 
\author{David P. Blecher}
\thanks{* Supported by a grant from the NSF}
\thanks{The contents of this paper were announced at
the January 1999 meeting of the American Mathematical Socety.}
\address{Department of Mathematics\\University of Houston\\Houston,
TX 77204-3476 }
\email{dblecher@@math.uh.edu}\maketitle
 
\maketitle
  
\vspace{10 mm}
 
\begin{abstract}
We continue our study of the general theory of
possibly nonselfadjoint algebras of operators on a Hilbert 
space, and modules over such algebras,
developing a little 
 more technology to connect `nonselfadjoint
operator algebra' with the C$^*-$algebraic framework.
More particularly, we make use of the {\em universal}, or
{\em maximal}, C$^*-$algebra generated by an operator algebra,
and C$^*-${\em dilations}.
This technology  is quite general,
however it was developed to solve some problems 
arising in the theory of Morita
equivalence of operator algebras, and as a result
most of the applications given here (and in 
a companion paper) are
to that subject.  Other applications given here are  
to extension problems for module maps, and
characterizations of C$^*-$algebras.
\end{abstract}
 
\pagebreak
\newpage
 
\setcounter{section}{0}

\section{Introduction - Modules over operator algebras}

In what follows $\A$ is a
possibly nonselfadjoint
operator algebra, that is, a
 general norm closed algebra of operators 
on a Hilbert space.  We shall assume that $\A$ has a
contractive approximate identity (c.a.i.).  Thus
any C$^*-$algebra is an operator algebra.
The {\em general theory} of operator algebras, and 
of representations of, and modules over, such algebras, 
is lamentably sparse.
This is in contrast to the selfadjoint case, namely
the C$^*-$algebra theory, and the contrast is easily seen
in the lack of certain fundamental tools which are 
available in the selfadjoint case, such as 
von Neumann's double commutant theorem.  
This paper is the latest in a series
in which we study
the class of all operator algebras and 
their modules, using the recent perspectives
and techniques of `operator space' theory. 
One of the basic points of this latter theory (see \cite{E}), is 
that for many purposes
it is not sufficient to study linear spaces
of operators between Hilbert spaces in the classical functional
analytic  framework, namely in terms of norms and bounded linear 
maps.  One must use `matrix norms' and completely bounded linear
maps.  
We refer the reader to 
\cite{Ar,Pis,MK}, for background on operator spaces and 
operator algebras and a description of some 
other work in this area, and to 
\cite{BK} for a leisurely introduction and survey of
our work.   Our main purpose here is to expose some   
more connections between `nonselfadjoint
operator algebra' with the C$^*-$algebraic framework. 
Hitherto
many researchers seem to have assumed that there
is only one important C$^*-$algebra associated 
with a nonselfadjoint operator algebra $\A$, namely 
the C$^*-$envelope of $\A$.  In fact there is a
lattice of C$^*-$algebras generated by $\A$. 
The C$^*-$envelope, being the `smallest', is the easiest
 to concretely get one's hands on,
 and has many wonderful 
properties.  However, the maximal C$^*-$algebra 
C$^*_{max}(\A)$ generated by $\A$, which we concentrate
on here, has very useful 
properties which the C$^*$-envelope lacks, and 
is for some purposes more important.  

In this paper we study two kinds of 
representations of a nonselfadjoint algebra 
$\A$.  The first kind are the
completely contractive representations $\pi$ of $\A$ on a
Hilbert space $H$ say.  Then $H$ is naturally a left
$\A$-module: we shall refer to
such a module as  a {\em Hilbert} $\A$-{\em module}. 
Perhaps a better name might be {\em completely contractive 
Hilbert} $\A$-{\em module}, but we will use the shorter 
name since we do not care about
any other kind here.   
It is not assumed in this paper 
that such modules are nondegenerate\footnote{In this paper,
for a (left) Banach module $X$ over $\A$ we assume 
$\Vert a x \Vert 
\leq \Vert a \Vert \Vert x \Vert$ for $a \in \A, x \in X$.  
A left Banach module is said 
to be {\em essential} or {\em nondegenerate} if
$\{\sum_{k=1}^n a_k x_k : n \in \N , \; a_k \in \A,
\; x_k \in X \}$ is
dense in $X$.  This is equivalent to saying that for any c.a.i.
 $\{ e_\alpha \}$ in $\A$,
$e_\alpha x \rightarrow x$ for all $x \in X$.
Banach modules are not assumed nondegenerate here
unless explicitly stated.}, 
unless we explicitly say so.
If $\A$ is a C$^*-$algebra, it is folklore
(but also follows from our Theorem \ref{im1}) that
contractive representations on Hilbert space are
$*-$representations, and thus automatically
completely contractive.

The  second type of representation of $\A$, which is more 
general than the first type, corresponds to what is known
as an {\em operator $\A$-module}, and is explained in more 
detail below.   
We shall explore the 
connections between  
the study of
Hilbert and operator modules over $\A$, and those over 
$\C$,
where $\C = C^*_{max}(\A)$ is the maximal, or
 universal, C$^*-$algebra 
generated by $\A$.  We reserve the symbol 
$\C$ for this C$^*-$algebra throughout.
In  \S 2 we show how to construct $\C$ and give
some examples.  
It turns out, although this is not as
obvious as at first glance it appears to be,
that the class of operator
modules over $\C$, is a subcategory
of the class of operator modules over $\A$.  We derive
this in section 3
from some general (but apparently new)
facts about Banach modules over 
C$^*$-algebras.  Moreover 
every Hilbert or operator $\A$-module `dilates'
to an operator module over $\C$.  This process
is studied in \S 3, the central section of this paper,
which shows how the category of modules over
a nonselfadjoint operator algebra $\A$, and the
category of modules over
any C$^*-$algebra generated by $\A$, are related.
As a first application of this and some related ideas,
 we give in \S 4, a characterization of C$^*-$algebras 
amongst the operator algebras, in terms of injectivity of
certain modules, and in terms of the above dilations.

Of course the main motivation for the machinery 
developed here, is that certain problems concerning 
nonselfadjoint algebras
should be solvable by transferring them to the selfadjoint
framework, and then using C$^*-$algebra techniques.
An example of this principle is given in 
a companion paper \cite{BlFTM2}, where we use all the results
developed here in \S 3, to generalize the
main result of \cite{BlFTM1} to nonselfadjoint 
operator algebras.  This completes the circle of ideas begun in
\cite{BMP} concerning strong Morita equivalence of 
operator algebras.    We devote  \S 5 of the present paper
to various other
connections between C$^*-$dilations and strong
and `weak' Morita equivalence of operator algebras. 
The reader may find this a rather complicated application of 
the dilation, however it was our motivation for developing
the technology of the earlier sections.  We have no doubt 
that other, more simple, applications of this technology 
will follow in the course of time.  At present we are 
working on some connections between these ideas and
some interesting problems concerning function algebras
\cite{BJ}.

In \S 5, study of the C$^*-$dilation leads us to
define a new notion of Morita equivalence of operator
algebras, which we call `strong subequivalence',
which has many of the features one associates with
strong Morita equivalence.   It is called
strong subequivalence because, basically, it is an equivalence
which may be dilated to a strong Morita equivalence of the
generated C$^*$-algebras.
Strong
Morita equivalence implies strong subequivalence, but
the converse is false.  Strong subequivalence is thus a 
weaker notion 
than strong Morita equivalence, and thus is easier to 
check in particular examples, while having many of the same
consequences.  However, we show that strong
Morita equivalence of operator algebras,
is the same as strong subequivalence
when the last-mentioned dilation is to the
{\em  maximal} generated C$^*-$algebras.   This may be
viewed as a new characterization of  strong
Morita equivalence of operator algebras. 
We also show that a subcontext of a C$^*-$algebraic 
strong Morita equivalence is dilatable if and only if
it preserves the C$^*-$algebraic
weak Morita equivalence.

In \S 6 we study a class of operator 
modules, and C$^*$-modules, which can be associated
with any operator space or operator module.  We also define,
using the maximal C$^*-$dilation,
 a canonical operator algebra $\U(X)$, which we call the upper
linking algebra, associated 
with any operator bimodule $X$, which has an appropriate
universal property for completely contractive bimodule 
maps defined on $X$.

Let us begin 
by establishing the common symbols and notations.  
We shall use 
operator spaces quite extensively, and their connections to 
operator modules.  We refer the reader to 
\cite{BMP}, \cite{Bna} and \cite{BK}
for missing background.

Suppose that $\pi$ is a 
completely contractive representation of $\A$ on 
Hilbert space $H$, and that $X$ is a closed subspace of 
$B(H)$ such that $\pi(\A) X \subset X$.  Then $X$ is a left $\A$-module.  
We say that such $X$, considered as
an abstract operator space and a left $\A$-module, is a
 left {\em operator module} over $\A$.  
By considering $X$ as an abstract operator 
space and module, we may forget about the particular $H, \pi$ 
used.
We shall assume in future, unless we explicitly say
to the contrary, that the module action on an operator 
module $X$ is nondegenerate.
It is sometimes useful, 
and equivalent, to allow $X$ in the definition
above, to be a subspace of 
$B(K,H)$, for a second Hilbert space $K$.
The advantage of this is that it will allow  
$H  = [XK]^{\bar{}}$ if we wish.  (The notation
$[YZ]^{\bar{}}$ in this paper will mean the closure
of the linear span of products of terms in $Y$ and $Z$).
An obvious modification of a
 theorem of Christensen-Effros-Sinclair \cite{CES} tells us that
the operator modules are exactly the operator spaces which are 
(nondegenerate) left $\A$-modules, such
that the module action 
satisfies $\Vert a x \Vert \leq \Vert a \Vert \Vert x \Vert$
just as for a Banach module, except that now $a$ and $x$ may be 
square matrices of the same finite size, with entries in 
$\A$ and $X$ respectively.  In other words, the module action
is a `completely contractive' bilinear map (or equivalently,
the module action linearizes to a 
complete contraction $\A \otimes_h X \rightarrow X$,
where $\otimes_h$ is the Haagerup tensor product).  
We write $_{\A}OMOD$ for the category 
of left $\A$-operator modules.  The morphisms 
are $_{\A}CB(X,W)$, the {\em completely bounded} left 
$\A-$module maps.  Unless specified otherwise, when
$X, W$ are operator modules or bimodules, when we say
`$X \cong W$' , or `$X \cong W$ as operator modules', we 
mean that the implicit isomorphism is a completely isometric
module map.
If $X,W$ are left $\A$-operator modules then
$_{\A}CB(X,W)$ is an operator space, whose operator 
structure is specified by the natural (algebraic)
identification
$M_n(_{\A}CB(X,W)) \cong \; _{\A}CB(X,M_n(W))$. 

We let $_{\A}HMOD$ be the category   
of nondegenerate Hilbert $\A$-modules.
In \cite{BMP}
we showed how $_{\A}HMOD$ may be viewed as a subcategory of
$_{\A}OMOD$ (see the 
discussion at the end of Chapter 2, and after 
Proposition 3.8, there).
Briefly, if $H \in $ $_{\A}HMOD$,
 then if $H$ is equipped with its
Hilbert column operator space structure $H^c$, then 
$H^c \in $ $_{\A}OMOD$.
Conversely, if $V \in \; _{\A}OMOD$ is also a
 Hilbert
column space, then the associated representation $\A \rightarrow
B(V)$ is completely contractive and nondegenerate.
It is well known that for a linear map $T : H \rightarrow K$
between Hilbert spaces, the
usual norm equals the completely bounded norm of $T$
as a map $H^c \rightarrow K^c$.
Thus we see that the assignment 
$H \mapsto H^c$  embeds $_{\A}HMOD$ as a (full)
subcategory of 
$_{\A}OMOD$.  In future if a Hilbert space is referred to
as an operator space, it will be with respect to its 
column operator space structure, unless specified to
the contrary.

In \cite{BK} Lemma 8.1
we showed that if $\A$ is an operator algebra
with contractive approximate identity $\{ e_\alpha \}$, and
if $\D$ is any C$^*-$algebra generated\footnote{That is, 
$\D$ is a C$^*-$algebra generated by a completely isometric,
homomorphic, copy of $\A$.} by 
$\A$, then 
$\{ e_\alpha \}$ is a contractive approximate identity
for $\D$.  This fact will be used frequently.  In particular
if follows that the obvious action of $\A$ on $\D$ is
nondegenerate, so that $\D \in \; _{\A}OMOD$.

We usually choose work with left
modules here.  The right module versions, or bimodule 
versions, are mostly similar.  There is an important
principle which allows one to go  between right and 
left operator modules.  Namely, if $V$ is a left
module over $\A$, define $\bar{V} = \{ \bar{v} : v \in V \}$,
with the conjugate linear structure.  
Then $\bar{V}$ is a right module over $\A^*$.  Of 
course if $\A$ is a C$^*-$algebra, then 
$\A^* = \A$, otherwise one can view $\A^*$ as
the algebra of adjoints of $\A$ in any containing 
C$^*-$algebra.   There is
an obvious operator space structure to put on $\bar{V}$,
namely $\Vert [\bar{v}_{ij} ] \Vert_n =
\Vert [v_{ji}] \Vert_n$.   If $V$ is a left operator 
module over $\A$ then $\bar{V}$ is a right operator 
module, and we shall call it the conjugate operator
module.  

We end this section with a fairly obvious observation:

\begin{lemma}
\label{su}  Suppose that $\D$ is a C$^*$-algebra generated by
$\A$, that 
 $H$ and $K$ are Hilbert
$\D$-modules, and that $i : H \rightarrow K$
and $q : K \rightarrow H$ are contractive $\A-$module maps
with $ q \circ i = Id_H$.  Then $i$ and $q$ are $\D-$module
maps.  In particular, a unitary $\A-$module map
$u : H \rightarrow K$ is a
$\D-$module map.   \end{lemma}
 
For completeness we
give the easy proof.   By a basic fact about contractions
on a Hilbert space, we have $q = i^*$.  Let us write
$\rho$ and $\sigma$ for the representations of $\D$ on
$H$ and $K$ respectively.  Then
for $a \in \A , \zeta \in H , \eta \in K$ we have
$$\langle i(\rho(a)^* \zeta) , \eta \rangle =
\langle \zeta , \rho(a) q(\eta) \rangle =
\langle \zeta , q(\sigma(a) \eta) \rangle = \langle
\sigma(a)^* i(\zeta) , \eta \rangle \; .
$$
This shows that $i$ is a $\D-$module map. Similarly
$q$ is a $\D-$module map.

\section{The maximal C$^*-$algebra.}
   
In \cite{BP2} we defined the universal or
maximal C$^*-$algebra of an operator 
algebra $\A$, and it appeared again in \cite{BMN}.
Since it did not play a particularly significant role in those 
papers, we did not give a careful development.  We begin by
remedying this omission.  

\begin{definition}
\label{cmax}  If $\A$ is an operator algebra with contractive approximate
identity, then there exists a C$^*-$algebra $\C$ and a completely 
isometric homomorphism $i : \A \rightarrow \C$ such that $i(\A)$
generates $\C$ as a C$^*-$algebra, and such that if $\phi : A
\rightarrow \D$ is any completely contractive homomorphism into a
C$^*-$algebra $\D$, then there exists a (necessarily unique)
*-homomorphism $\tilde{\phi} : \C \rightarrow \D$ such that
$\tilde{\phi} \circ i = \phi$.  The C$^*-$algebra $\C$ is called
the maximal C$^*-$algebra generated by $\A$, and is 
written as $C^*_{max}(\A)$.
\end{definition} 
    
The existence and uniqueness
of such a universal object $(\C,i)$ is not difficult,
but since it is not written anywhere in the literature as far
as we are aware, we give the details.  
We may suppose that $\A$ has an identity of norm 1
(otherwise adjoin an identity in the usual way,
and let $\C$ be the C$^*$-subalgebra of
$C^*_{max}(\A^{1})$ generated by $\A$).  
Let $\E$ be the algebraic free product of $\A$ and $\A^*$,
which is clearly a $*$-algebra.  
We now use some basic facts from \cite{Ar} or
\cite{P} about completely
positive maps.  We recall that the operator algebra
$\A^{*}$, and indeed the operator system
$\A + \A^{*}$ does not depend on 
any particular Hilbert space that $\A$ is represented on.
Let $\theta : \A \rightarrow
\D$ be a  c.c. homomorphism 
into a C$^*-$algebra $\D$.  Let $\D'$ be the 
C$^*-$algebra generated by the range of $\theta$. 
Then $\theta$ extends to a completely
positive unital map $ \A + \A^* \rightarrow \D'$, 
which when restricted to $\A^*$ is a c.c. homomorphism
$\theta^*$.  Then $\theta \star \theta^* : \E \rightarrow
\D'$ is a $*$-representation.  In the usual way, $\E$ gives 
rise to a C$^*-$algebra $\C$ by taking the supremum
over all such $*$-representations.  Clearly $\A$ is
unitally completely isometrically embedded as a subalgebra
of $\C$, $\A$ generates $\C$, and $\C$ has the required 
universal property.  This gives the existence of $\C$.  
However,   $\C$ is clearly unique in the sense that
if $(\C',i')$ is any other pair with the 
property described in \ref{cmax} then there exists a unique
*-isomorphism $\pi : C \rightarrow C'$ with $\pi \circ i = i'$.

\vspace{5 mm}

\begin{proposition}
\label{Rem}
We have
 C$^*_{max}(\A_1 \star \A_2)
\cong C^*_{max}(\A_1) \star C^*_{max}(\A_2)$,
for operator algebras $\A_1$ and $\A_2$ with c.a.i.,
 where $\star$ is the operator algebra
free product of \cite{BP2}.
\end{proposition}

This follows immediately from the universal properties.
The analogous result for the maximal operator algebra
tensor product of \cite{PP} is certainly false, as may
be seen for example from \ref{t2} below, and 2.6 in 
\cite{PP}.

\vspace{5 mm}

\noindent {\bf Remark.}   For those who are 
familiar with  operator space theory, it
is tempting to think of $\C = C^*_{max}(\A)$ as an infinite
Haagerup tensor product of copies of $\A$ and $\A^*$.
Indeed it is tempting to think of elements in 
$\C$ , in the spirit of \cite{BP2}, as products $A_1 B_1
A_2 B_2 \cdots $ of
matrices, with $A_i$ from $\A$ and $B_i$ from $\A^*$.
From
this perspective one might be led to conjecture that
$\C \otimes_{h} \C \cong \C$ or
$\C \otimes_{h\A} \C \cong \C$.  The first conjecture is
false unless $\A = {\mathbb C}$
by \cite{B1} Theorem 1.  The second is also false, as we shall 
see later in \ref{cod}.

\vspace{5 mm}

Notice that in \ref{cmax} one may take w.l.o.g. the $\D$ there
to be $B(H)$ for a Hilbert space $H$.    
That is, we may take the $\phi$ in \ref{cmax} to be a completely 
contractive representation.  From this we see immediately that
Hilbert $\A$-modules are automatically Hilbert $\C$-modules
and vice versa.  Thus as 
{\em objects} $_{\A}HMOD = $ $_{\C}HMOD$.  However the morphisms
in these two categories
are not the same, since $\A$-intertwiners are not necessarily
$\C$-interwiners.  In fact it is clear that $_{\C}HMOD$ is a
subcategory of $_{\A}HMOD$.  

We remark that it seems interesting to
 transfer the language of the representation theory of
C$^*-$algebras to operator algebras.  Thus for example
we say that $\A$ (or a representation
$\phi$  of $\A$), is type I or CCR, and so on,
 if and only if $\C$ (or $\tilde{\phi}$) has this property.
For example, by results in \cite{FS}, the disk algebra is
NGCR.  This example is discussed further in \ref{AD}.  

In the rest of this paper we will take $H_u$ to be the
Hilbert space of the universal representation of $\C$, 
and will refer to $H_u$ as the universal representation 
of $\A$.  It is clear from C$^*-$algebraic representation
theory, that any nondegenerate
Hilbert $\A$-module is (completely)
isometrically $\A-$module isomorphic to a
complemented $\A$-submodule of a direct sum of 
copies of $H_u$.

\begin{example}
\label{AD}  
\end{example}
Consider $\A = A({\mathbb D})$, the disk algebra.
In this case 
$\C = C^*_{max}(\A)$ is the universal C$^*-$algebra 
generated by a contraction, which has been studied by many
researchers.  This is a noncommutative C$^*-$algebra,
generated by a non-normal contraction $z$, say.  

We found this example helpful in disposing of several
incorrect guesses we had concerning $C^*_{max}$.  
For example, one can use it to show that
if $S \in \; _{\A}CB(\C,\C)$, and
$S(a) = 0$ for all $a \in \A$, then $S$ is not necessarily
the zero map.
Define $L(c) = z c$
for $c \in \C$.  Clearly $L \in $$_{\A}CB(\C)$, but $L$ is not in
$_{\C}CB(\C)$ since $z$ is not normal.  
If we put $S = L$, notice that $S$ restricted to
$\A$ equals $r_z$, i.e. right multiplication by 
$z$.  Hence $S - r_z$ is a left $\A$-module map
on $\C$, is zero on $\A$, but is not the 
zero map. 

\vspace{5 mm}

\begin{example}
\label{t2}
\end{example}
Consider $\A = {\mathcal T}(2)$,
 the upper triangular 
$2 \times 2$ matrices.
Let $\A_0$ be its subalgebra consisting of those matrices
with repetition on the main
diagonal.   Then $C^*_{max}(\A_0)$ is the 
well known 
universal C$^*-$algebra generated by a nilpotent operator,
and $\A_0$ is the universal operator algebra 
generated by a nilpotent operator.
In \cite{Lor}, $C^*_{max}(\A_0)$ is shown to be
$\{ f \in M_2(C([0,1])) : f(0) \in {\mathbb C} I \}$,
also known as the cone over $M_2$.
In fact 
$C^*_{max}(\A) = \{ f \in M_2(C([0,1])) : 
f(0)$ is a diagonal matrix$\}.$  
We will prove this (and a little bit more).
It is convenient to work with 
the dense subalgebra $\E$ of the
 last C$^*-$algebra
consisting of matrices of the form 
$$
\left[ \begin{array}{ccl}
b_1 & b_2 \sqrt{t} \\
b_3 \sqrt{t}  & b_4
\end{array}
\right] \; \; 
$$
here `$t$' is the basic monomial on $[0,1]$,
and $b_i \in C([0,1])$.
Note that this last C$^*-$algebra
 is generated by the subalgebra
consisting of matrices
$$
\left[ \begin{array}{ccl}
\lambda_1  & \mu \sqrt{t} \\
0 & \lambda_2 
\end{array}
\right] \; \;
$$
where $\lambda_1,\lambda_2, \mu \in {\mathbb C}$.
This subalgebra is easily seen to be completely isometrically
isomorphic to ${\mathcal T}(2)$, and we will 
henceforth take $\A$ to be
this subalgebra.
Notice that if $T : K \rightarrow H$ is any
contractive operator between Hilbert spaces, then the
subspace of $B(H \oplus K)$ consisting of matrices
$$
\left[ \begin{array}{ccl}
\lambda_1 I & \mu T \\
0 & \lambda_2 I 
\end{array}
\right] \; \;
$$
with $\lambda_1,\lambda_2, \mu$ scalar, is an algebra.
The $C^*-$algebra it generates consists of all the
matrices of the form 
$$
\left[ \begin{array}{ccl}
p_1(T T^*) & p_2(T T^*) T \\
T^* p_3(T T^*)  & p_4(T^*T)
\end{array}
\right] \; \; ,
$$
where the $p_i \in C([0,1])$.  
It is easily checked that
the map
$$
\left[ \begin{array}{ccl}
p_1 & p_2 \sqrt{t} \\
p_3 \sqrt{t} & p_4
\end{array}
\right] \; \; \mapsto
\left[ \begin{array}{ccl}
p_1(T T^*) & p_2(T T^*) T \\
T^* p_3(T T^*)  & p_4(T^*T)
\end{array}
\right] \; \; ,
$$
is a *-homomorphism from $\E$ into 
$B(H \oplus K)$.  However it is clearly continuous - notice
for example, that 
$\Vert p_2(T T^*) T \Vert
= \Vert p_2(T T^*) (T T^*)^{\frac{1}{2}} \Vert
\leq \Vert p_2 \sqrt{t} \Vert_{[0,1]}$.
Hence it extends to a *-homomorphism on the 
containing C$^*-$algebra, and is consequently
completely
contractive.  By restriction we obtain a 
completely contractive homomorphism from 
$\A$ into $B(H \oplus K)$.  Conversely,
any nilpotent operator on a Hilbert space $L$,
or any nondegenerate
contractive representation of $\A$,
immediately gives a decomposition of $L$ as
$H \oplus K$, and an operator $T$ as above,
with respect to which we are
again in the above situation.

The above shows that $\{ f \in M_2(C([0,1])) :
f(0)$ is a diagonal matrix$\}$ may be characterized as the
universal unital 
C$^*-$algebra generated by two contractions
$x,v$ with relations $x^2 = 0, 
v^2 = v , v x = x , x v = 0$.  This fact is no doubt
well known.  
This seems related to  2.6 in \cite{PP}
which says, loosely, 
that a commutant lifting theorem for general 
operator algebras follows from
knowing a certain result for ${\mathcal T}(2)$.

We generalize the previous example in the final section 
of our paper.

\section{Operator modules over a generated C$^*-$algebra
and C$^*-$dilations.}

This is the central section of this paper, in which
we show how a category of modules over
a nonselfadjoint operator algebra $\A$, and the
category of modules over
a C$^*-$algebra generated by $\A$, are related 
by an interesting pair of adjoint functors.  All 
the results developed here are heavily relied on
in \cite{BlFTM2}, and later in the present paper,
and should be useful in many 
other situations.    

We begin this section
 with some general facts about Banach modules over
C$^*-$algebras which, as far as we are
aware, are new.  In \cite{B1} we proved
the following result:

\begin{theorem}
\label{im1}  Let $\D$ be a C$^*-$algebra, $\B$ a Banach algebra,
and $\theta : \D
\rightarrow \B$ a contractive homomorphism.  Then the
range of $\theta$ is norm-closed, has a contractive
approximate identity, and possesses an involution with respect
to which it is a C$^*-$algebra.
\end{theorem}

Thus if $V$ is a left Banach module over a C$^*-$algebra $\D$,
and if we let $\theta : \D
\rightarrow B(V)$ be the associated contractive homomorphism
then the range of $\theta$ is a C$^*-$algebra.

\begin{theorem} \label{unq}  
Suppose that $V$ is a Banach module
over
an operator algebra $\A$ with contractive approximate
identity.  Write $\theta : \A \rightarrow B(V)$
for the associated  homomorphism.  Suppose that
$\D$ is any C$^*-$algebra generated by $\A$.
Clearly the
$\A$-action on $V$ can be extended to
a $\D$-action with respect to which $V$ is a Banach
$\D$-module if and only if
$\theta$ is the restriction to $\A$ of
a contractive homomorphism $\phi : \D \rightarrow B(V)$.
This extended $\D$-action, or equivalently the
homomorphism $\phi$, is unique if it exists.
\end{theorem}

\begin{proof}
Only the uniqueness requires proof.
We shall require some facts about Banach algebras
which may be found in \cite{BD}.
Suppose that $\phi_1$ and $\phi_2$ are two contractive
homomorphisms $\D \rightarrow B(V)$, extending $\theta$.
By Theorem \ref{im1},
the ranges $\E_1$ and $\E_2$ of $\phi_1$ and $\phi_2$ are
each C$^*-$algebras, but with possibly different involutions.
We will write these involutions as $*$ and $\#$
respectively.  With respect to these involutions
$\phi_1$ and $\phi_2$ are `*-homomorphisms'.
Choose a c.a.i. $\{ e_\alpha \}$ for $\A$,
and let $\B = \{ T \in B(V) : T \theta(e_\alpha) \rightarrow
T$, and $\theta(e_\alpha) T \rightarrow T \}$.   Then
$\B$ is a Banach algebra with c.a.i. $\{ \theta(e_\alpha) \}$.
If $\F$ is a Banach algebra with c.a.i. we define
 $\F^{1} = \F$ if $\F$ is unital, otherwise we
let it be the unitization of $\F$, with its `multiplier
norm'.   In the nonunital case, it is easy to see that
$\F^{1}$ may be defined equivalently to be the
subalgebra of $\F^{**}$ generated by $\F$ and a
weak*-limit point of the c.a.i.  In any case,
the `unitized' C$^*-$algebras 
$\E_1^1$ and $\E_2^1$ 
may be viewed as subalgebras of $\B^1$,
with the same unit.
Let $a \in \A$, and let $f$ be a state on $\B$
(or equivalently on $\B^{1}$).  Then for $k = 1,2$,
$f$ restricted to $\E_k$ is a state on $\E_k$.  Thus
$f(\phi_1(a)^*) = \overline{f(\phi_1(a))} =
\overline{f(\phi_2(a))} = f(\phi_2(a)^{\#})$.
Thus $u = \phi_1(a)^* - \phi_2(a)^{\#}$
is a Hermitian element in $\B$ (or $\B^{1}$)
with numerical radius $0$, and consequently $u = 0$.
Therefore $\phi_1 = \phi_2$ on $\D$.
\end{proof}

From this we obtain the following `rigidity' result:

\begin{corollary}
\label{cor2}  Let $\D$ be a C$^*-$algebra generated by
an operator algebra $\A$.  If $V_1$ and $V_2$ are two
Banach $\D$-modules, and if $T : V_1 \rightarrow V_2$ is an
isometric and surjective $\A$-module map,  then
$T$ is a $\D$-module map.
\end{corollary}

\begin{corollary}  Let $\D$ be a C$^*-$algebra generated by
an operator algebra $\A$.  The category of Banach modules over
$\D$ is a subcategory of the category of Banach modules over
$\A$.  Similarly, $_{\D}OMOD$ is a subcategory of
$_{\A}OMOD$, and $_{\D}HMOD$ is a subcategory of
$_{\A}HMOD$.
\end{corollary}

Thus the `forgetful
functor' from the category of Banach (or operator, or 
Hilbert)
modules over $\D$, to the same category over $\A$, is
unambiguous (i.e. one-to-one), and embeds the
first category as a subcategory of the second.
In more flowery language \cite{Cat} it turns out that
the subcategory is `reflective'.
We regard it as one of the significant
open problems in this area
to find a good test for when an $\A$-operator module $V$
possesses an extended $\C$-module action.

In the remainder of this section
we discuss the `$\D$-dilation' of an $\A$-operator 
module $V$, where $\D$ is a C$^*$-algebra generated by $\A$.
We shall see that, in 
the language of category theory, the $\D$-dilation is
the left adjoint of the aforementioned forgetful functor  
from $_{\D}OMOD$ to $_{\A}OMOD$.
In fact it is a simple 'change of rings'.
The word 'dilation', and some of its useful properties,
 was first introduced in work of  Muhly and
Na \cite{MN,Na}, in 
the case when $\D$ is the C$^*$-envelope of $\A$.  We
will indicate as we go along, any 
overlap with their work.  The dilation was also used, but 
not explicitly named as such, in   
\cite{Bhmo} and \cite{BMN}.  

\begin{definition}
\label{upp} 
A pair $(E,i)$
is said to be
a $\D$-{\em dilation} of a left
$\A$-operator
module $V$, if both of the following hold:
\begin{itemize}
\item [(*)] $E$ is a left $\D$-operator
module and   
$i : V \rightarrow E$ is a 
completely contractive $\A$-module map, 
\item [(**)]  
For any left $\D$-operator module $V'$, and any
completely bounded
$\A-$module map $T : V \rightarrow V'$, there
exists a 
unique completely bounded
$\D-$module map $\tilde{T} :
E \rightarrow V'$ such that $\tilde{T} \circ i = T$, and
also
$\Vert \tilde{T} \Vert_{cb} = \Vert T \Vert_{cb}$.
\end{itemize}
\end{definition}

This is a universal property in the sense that  if
$(E',i')$ are any pair
satisfying (*) and (**),
then there
exists a unique completely isometric $\D$-module isomorphism
$\rho : E \rightarrow E'$ such that $\rho \circ i = i'$.
We will postpone the existence of the $\D$-dilation to the
next lemma.

The `uniqueness' assertion in (**) is equivalent to
saying that  $i(V)$ generates $E$ as a $\D-$operator
module in the obvious sense (namely,
that there are no nontrivial closed $\D$-submodules
of $E$ which contain $V$).  To see this let  
$E' = [\D i(V)]^{\bar{}}$, and consider the quotient
map $Q : E \rightarrow \frac{E}{E'}$.

The $\D$-dilation $(E,i)$ is clearly
the unique pair satisfying (*),
such that for all $\D$-operator modules $V'$,
the canonical map
$i^* :
\; _{\D}CB(E,V') \rightarrow $$_{\A}CB(V,V')$, given by
composition with $i$, is an
isometric 
isomorphism.  Since $M_n(CB(X,Y)) = CB(X,M_n(Y))$ for operator
spaces, it follows that $i^*$
being an
isometry for all such $V'$ implies that it is a
complete isometry.
Thus the $\D$-dilation $E$ of $V$ satisfies: 
$$_{\D}CB(E,V') \cong \; _{\A}CB(V,V') \; \; .  \eqno{(***)} $$
completely isometrically.  In the case that 
$\D$ is the C$^*$-envelope of $\A$, a part
of this assertion was observed by Muhly and Na.
In fact, what
this result says in the language of elementary category 
theory \cite{Cat}, is that the $\D$-dilation is the left
adjoint of the forgetful functor from $_{\D}OMOD$ to
$_{\A}OMOD$ (discussed at the end of \S 1).  Of course, 
either of the two  
compositions of this forgetful functor and the 
$\D$-dilation
is not the identity.
Another good name for what we call the $\D-$dilation
might be the `$\D$-adjunct'.  In flowery language,
this {\em adjunction} makes $_{\D}OMOD$ a {\em
reflective subcategory}
of $_{\A}OMOD$.

The following shows that we may 
take $E$ to be the
Haagerup module tensor $\D \otimes_{h\A} V$ .
See \cite{BMP} for the definition of the module Haagerup
tensor product $\otimes_{h\A}$, as well as for
its basic properties, such as the fact that it is
associative, functorial, and that $\A \otimes_{h\A} V \cong V$.
We note that
since $\A \otimes_{h\A} V \cong V$, there is a canonical
completely contractive $\A$-module map
$i : V \rightarrow \D \otimes_{h\A} V$ .

\begin{lemma}
\label{up}  For any left operator module $V$ over $\A$, 
the $\D$-operator module $E =  \D \otimes_{h\A} V$
is the $\D$-dilation of $V$.  
\end{lemma} 

\begin{proof}  If  $T : V \rightarrow V'$ is 
as above, then by the functoriality of the Haagerup 
tensor product, $Id_{\D}
\otimes T : \D \otimes_{h\A} V 
\rightarrow \D \otimes_{h\A} V'$ is
completely bounded.  Composing this with the module action 
$\D \otimes_{h\A} V' 
\rightarrow V'$ gives the required map $\tilde{T}$.
Its easy to see that $\tilde{T}$ has the right properties.   
The uniqueness assertion is obvious.  
\end{proof}

\vspace{4 mm} 

We now make some observations which will be 
important to us.  First, notice that
it is not necessary 
that the $V'$ be nondegenerate in (**) above,
since one may always
replace $V'$ with its `$\D$-essential submodule'.  Note that any 
$T$ as above maps into this essential submodule of $V'$.   
Secondly, observe that by the 
Christensen-Effros-Sinclair result, 
it suffices to take 
$V' = B(H,K)$ in (**),
where $K$ is a Hilbert 
$\D$-module and $H$ is a Hilbert space.
The next theorem shows that with a natural 
qualification, one may as well 
take $V' = K$.  

It was probably first noted by Effros
that $_{\D}CB(F,B(H,K))$ is a dual operator
space, if $F$ is a left $\D$-operator
module.  Using basic results about
operator spaces, it can be shown that its operator space
predual may be written
as  $K^r \otimes_{h\D} F \otimes_h H^c$, where
$K^r$ is the operator dual of $K$.  The duality pairing
is the obvious one, namely,
$\langle T , \psi \otimes x \otimes \zeta \rangle
= \langle T(x) (\eta) , \psi \rangle$, for 
$T \in \; _{\D}CB(F,B(H,K)), x \in X, \zeta \in H ,
\psi \in K^*$. 
Similarly
for $_{\A}CB(V,B(H,K))$.  Note that there is
a canonical complete contraction
$S : K^r \otimes_{h\A} V \rightarrow K^r \otimes_{h\D} E$ formed
from the composition of the following maps:
$$K^r \otimes_{h\A} V
\overset{Id \otimes i}{\rightarrow} K^r \otimes_{h\A} E
\cong K^r \otimes_{h\D} \D \otimes_{h\A} E 
\rightarrow K^r \otimes_{h\D} E \; \; .$$
The last map in this sequence comes from the multiplication
$\D \times E \rightarrow E$.
We then get a map
$S_1  = 
S \otimes Id_H : K^r \otimes_{h\A} V \otimes_h H^c \rightarrow
K^r \otimes_{h\D} E \otimes_h H^c$.
It is easy to check that $S_1^*$ is what we called $i^*$
earlier.  Hence $i^*$
is an isometric isomorphism if and only if $S_1$ is
an isometric isomorphism.  

If  $T : X \rightarrow Y$
is a contraction (resp. isometry) between operator spaces,
then we will say that $T$ is a 
{\em row contraction} (resp. {\em row isometry}),
if $\Vert  [T(x_1) \; T(x_2) \cdots T(X_n)] \Vert \leq 
$ (resp. $=$)$\Vert  [x_1 \cdots x_n] \Vert$ for all $n$
and $x_i \in X$.  The following is mostly
in \cite{Ma}, but for completeness we give a proof.

\begin{lemma}
\label{co}
Let $T : X \rightarrow Y$ be a linear map between 
operator spaces.  The following are equivalent:
\begin{itemize}
\item [(i)]  $T$ is a row contraction.
\item [(ii)] $T^*$ is a row contraction.
\item [(iii)]  For all Hilbert spaces $H$,
$T \otimes I_H : X \otimes_h H^c \rightarrow
Y \otimes_h H^c$ is a contraction.
\end{itemize}
\end{lemma}

\begin{proof}  
Note that by definition of $\otimes_h$, (i) implies
(iii) (and in fact one may replace $H^c$ by any
operator space).  Put $H = C_n$ in (iii), and 
observe that $(X \otimes_h C_n)^* \cong CB(X,R_n) \cong
R_n(X^*)$.  Dualizing $T \otimes I_n$ now yields
(ii).  Since, therefore, (i) implies (iii), we see 
that (ii) implies that $T^{**}$ is a row contraction,
so that $T$ is also.
\end{proof}

Putting the observations above together,
we obtain the equivalence of (**) and condition
(i) or (ii) or (iii) below:

\begin{theorem}
\label{eno}  Suppose a pair $(E,i)$ satisfies (*).
Then $(E,i)$ satisfies (**), and consequently is the
$\D$-dilation of $V$, if and only if 
one of the following properties holds: 
\begin{itemize}
\item [(i)] the canonical  
map $i^* : \; _{\D}CB(E,K) \rightarrow $$_{\A}CB(V,K)$ 
defined above
is a completely 
isometric isomorphism, for all Hilbert $\D$-modules $K$. 
\item [(ii)] the canonical 
map $i^* : \; _{\D}CB(E,K) \rightarrow $$_{\A}CB(V,K)$
is a `row-isometric' isomorphism, 
for all Hilbert $\D$-modules $K$. 
\item [(iii)]  The canonical map
 $S : K^r \otimes_{h\A} V \rightarrow
K^r \otimes_{h\D} E$ 
defined above
is a `row-isometry',
for all nondegenerate Hilbert $\D$-modules $K$.
\end{itemize}
It is sufficient in (i) to take $K$ to be the universal 
representation of $\D$.
\end{theorem}
 
\begin{proof}
Only the last part still requires proof.
Every nondegenerate
Hilbert $\D$-module $K$ is a complemented 
submodule of a direct sum of $\gamma$ copies of the
universal representation, where $\gamma$ is some
cardinal.  Thus
 the last assertion of the theorem reduces to proving that:
if  the restriction map gives a complete
isometry $_{\D}CB(E,H) \cong \; _{\A}CB(V,H)$, then 
also $_{\D}CB(E,H^\gamma)
 \cong \; _{\A}CB(V,H^\gamma)$ completely 
isometrically.  One way to see this is to
first check that  
$_{\D}CB(E,H^\gamma) \cong M_{\gamma,1}(_{\D}CB(E,H))$
(see \cite{ERbimod}), and similarly
$_{\A}CB(V,H^\gamma) \cong M_{\gamma,1}(_{\A}CB(V,H))$.
Here $M_{\gamma,1}(X)$, for an operator space
$X$, is the collection of `columns' of length 
$\gamma$ with entries 
in $X$, whose truncated finite subcolumns are uniformly 
bounded.
\end{proof}

\vspace{4 mm}

The last statement of the previous theorem
is used in \cite{BlFTM2}.

If $V$ is an $\A-\B$-operator bimodule, where $\B$ is
a second operator algebra with c.a.i., then for any 
$\A$-operator module $V'$, the space $_{\A}CB(V,V')$ is
naturally a (not necessarily nondegenerate)
left $\B$-operator module with respect to
the action $(bT)(v) = T(vb)$.  We
define, as in \cite{BMP} \S 2, 
the space
$_{\A}CB^{ess}(V,V')$ to be the $\B$-essential subspace.
This equals $\{ T \in \; _{\A}CB(V,V') : 
T r_{f_\beta} \rightarrow T \}$, where $\{f_\beta\}$ is a
c.a.i. for $\B$, and $r_b$ is the operation of right multiplication
with an element in $\B$.  The maps in $_{\A}CB^{ess}$ we will 
refer to as `$\B$-essential'.
An important motivation for these 
spaces come from the theory of C$^*-$modules, where 
the `imprimitivity C$^*-$algebra'
or `algebra of `compact' operators', coincides with
$_{\A}CB^{ess}$.  We will need the next result later
and also in \cite{BlFTM2}:

\begin{theorem}
\label{bim}
For $V$ an $\A-\B$-operator bimodule, and for any 
Hilbert space $H$ and any (nondegenerate) Hilbert
$\A$-module $K$, we have that 
$_{\A}CB^{ess}(V,B(H,K))$ is weak*-dense in 
$_{\A}CB(V,B(H,K))$.
Moreover, if $\D$ is a C$^*$-algebra generated by $\A$, 
and if $(E,i)$ is a $\D-\B$-operator bimodule 
and a $\A-\B$-module map $i : V \rightarrow E$ whose 
range generates $E$ as a $\D$-operator module,
then $E$ is the $\D-$dilation of $V$ if and only if
$E$ satisfies (**) for $\B$-essential maps.
Moreover, the characterizations of the $\D$-dilation
in (i) and (ii) of the previous theorem, 
remain valid with $CB$ replaced by $CB^{ess}$.
\end{theorem}

\begin{proof}
If $T \in \; _{\A}CB(V,B(H,K))$, then the bounded 
net $\{ f_\beta T \}$ has a weak*-convergent subnet, 
which easily converges weak* to $T$.  That proves the first 
assertion.  Next notice that if $(E,i)$ satisfy (*), and if
$i^*$ is the canonical map
$_{\C}CB(E,B(H,K))
\rightarrow \; _{\A}CB(V,B(H,K))$
above, then $i^*(T)$ is $\B$-essential if and only if
$T$ is $\B$-essential.  From this it is easy to see the
`$\implies$' direction.  Conversely, given a
complete contraction $T \in 
\; _{\A}CB(V,B(H,K))$, then $f_\beta T$ lifts to a
complete contraction $S_\beta$ in $_{\C}CB(E,B(H,K))$.
A weak*-accumulation point of the $S_\beta$ will be
the desired extension of $T$.
We leave it to the reader to
fill in the remaining details.
\end{proof}

\vspace{4 mm}

\begin{lemma} \label{cov}
If $V$ is a left $\A$-operator module,
and if $\D$ is a C$^*-$algebra generated by $\A$,
then the following are equivalent:
\begin{itemize}
\item [(i)]  there exists
 a $\D$-operator module $V'$ and
a completely isometric $\A$-module map $j : V
\rightarrow V'$, and 
\item [(ii)]  the canonical $\A$-module map 
$i : V \rightarrow
 \D \otimes_{h\A} V$, is a complete isometry.
 \end{itemize}
\end{lemma}

\begin{proof}  Suppose that $m$ is the module action on
$V'$.  We have the following sequence of
canonical complete contractive $\A$-module maps:
 $$
V \overset{i}{\rightarrow} \D \otimes_{h\A} V 
\overset{Id \otimes j}{\rightarrow}
 \D \otimes_{h\A} V' \overset{m}{\rightarrow} V' \; \; .
$$
These maps compose to $j$, which yields the assertion.
\end{proof}

The idea of this last lemma was noticed by Muhly and Na
in the case that $\D$ is the C$^*-$envelope C$^*_e(\A)$
of $\A$.  We will refer
to the C$^*_e(\A)$-dilation as the `minimal C$^*$-dilation'.
In the case that $\D = \C = C^*_{max}(\A)$,
we call $\C \otimes_{h\A} V$ the `maximal C$^*$-dilation'.
A major reason for the usefulness of the latter
is the following, which follows immediately from
the previous result, the Christensen-Effros-Sinclair 
representation of operator modules, and the fact
that every Hilbert ${\A}$-module is a Hilbert
$\C$-module.

\begin{corollary}
\label{sub}  For any left $\A$-operator module $V$,
the  canonical $\A$-module map $i : V \rightarrow
 \C \otimes_{h\A} V$, is a complete isometry.
\end{corollary}

\vspace{5 mm}

We will regard $V$ henceforth 
as an $\A-$submodule of $\C \otimes_{h\A} V$.

There is
obviously an analogous C$^*-$dilation for right operator 
modules, or for operator bimodules.  The results in this 
section carry through without difficulty to these cases.

\section{Injectivity and characterizations of 
C$^*-$algebras.}

We now turn to some natural questions about injectivity,
C$^*-$dilations, and Hilbert modules which seem to be 
related.  Some of the results in this section may be known to 
experts, but it seems worthwhile to have them in print.

We will say that a (left) $\A$-operator module
$Z$ is (left) $\A$-injective if
whenever $V_2$ is a (left) $\A$-operator module
with closed submodule $V_1$, then every completely bounded
$\A$-module map $T : V_1 \rightarrow Z$ has a 
completely bounded
$\A$-module map extension
$\tilde{T} : V_2 \rightarrow Z$, 
with $\Vert T \Vert_{cb} = \Vert \tilde{T} 
\Vert_{cb}$.
Other authors do not require this 
last condition to hold,
and perhaps a better name for our property would be 
{\em 1-injective}.
Wittstock showed in \cite{Wi}
that if $\D$ is a unital C$^*-$subalgebra 
of $B(H)$ then $B(H)$ is $\D$-injective.   A rather
different proof may be found in \cite{Su}
(Suen uses bimodules, but the left module case 
can be easily obtained from his result by standard 
tricks).  The following consequence 
is fairly trivial, but we 
don't recall seeing it in the literature.  Another possible
proof of it,
using Suen's method, is described after Theorem \ref{pr}.   
We reaffirm that  we do not assume that 
Hilbert modules are nondegenerate, unless this is explicitly 
stated:

\begin{theorem}
\label{wit}  For any Hilbert module $H$ 
over a C$^*-$algebra $\D$,
$B(H)$ is (left) $\D$-injective.  More generally,
for any other Hilbert
space $N$, $B(N,H)$ is left $\D$-injective,
and $B(H,N)$ is right $\D$-injective.
\end{theorem}

\begin{proof}  By adjoining $I_H$ to $\D$, 
Wittstock's result fairly obviously extends to the
case when $\D$ is a nonunital 
C$^*-$subalgebra
of $B(H)$ acting nondegenerately on $H$.
Hence if
$H_u$ is the universal representation
of $\D$, and if
$K$ is a direct sum of copies of $H_u$, then 
  $B(K)$ is $\D$-injective.
However, every nondegenerate
Hilbert $\D$-module $H$ is a $\D$-complemented
submodule of such a $K$, and if $P$ is the 
$\D$-module projection onto $H$,
then $P B(K) P \cong B(H)$ as 
$\D$-operator modules.  Thus $B(H)$ is $\D$-injective.

If $H$ is not nondegenerate, we let $H'$ be the
essential part of $H$.  To show that $B(H)$ is injective,
 is sufficient to show that $B(H,H')$ is
$\D$-injective, since any $\D$-module map 
$T$ into $B(H)$ has range inside $B(H,H')$.
We may assume $H'$ is nontrivial,
otherwise the result is clear.   
However,  by a routine Hilbert space
cardinality argument
$B(H,H')$ may be regarded as a
$\D$-complemented submodule of $B(K,K)$ where $K$ is
a large enough direct sum of copies of $H'$.  

Finally, the $B(N,H)$ case is clear from the above,
whereas the right injectivity of $B(H,N)$ follows
from the left injectivity of $B(N,H)$ by noting 
that $B(N,H)$ is the `conjugate operator module' of
$B(H,N)$ .
 \end{proof}

The connection between injectivity and
dilations is explained by:  

\begin{proposition}
\label{buy}
Suppose that $V_2$ is an $\A$-operator module
with closed submodule $V_1$ .  Suppose that $\D$ is a
C$^*-$algebra generated by $\A$.  Then the following
are equivalent:
\begin{itemize}
\item [(i)] The  canonical map 
from the $\D$-dilation of $V_1$ to
$\D$-dilation of $V_2$
is a complete isometry.
\item [(ii)]  For every $\D$-injective module $\B$, and every 
completely bounded
$\A$-module map $T : V_1 \rightarrow \B$, then $T$ has a
completely bounded
$\A$-module map extension
$\tilde{T} : V_2 \rightarrow \B$,
with $\Vert T \Vert_{cb} = \Vert \tilde{T} \Vert_{cb}$.
\item [(iii)]  
For every Hilbert $\D$-module $K$,
the canonical map $K_r \otimes_{h \A} V_1 
\rightarrow K_r \otimes_{h \A} V_2$ is a complete isometry,
where $K_r$ is the operator dual of $K$).  
\item [(iv)]  Same as (iii), but with a single Hilbert
module, namely the Hilbert
space of the universal representation of $\D$.
\end{itemize}
\end{proposition}

\begin{proof}  Note that just as in the Remarks after
\ref{up}, it suffices to take $\B$ in (ii) to be 
$B(H,K)$, where $H$ is a Hilbert space, and $K$ is an
Hilbert $\D$-module.  By an argument similar to that
given in those same Remarks (the main difference being
that the map $i^*$ there is a complete quotient map), 
this is equivalent to (iii) (in fact, one may 
replace the word `complete' in (iii) with `row').
To see that (iv) implies (iii), we first observe that as in
\ref{wit}, we may assume $K$ is nondegenerate.  Using the 
functoriality of $\otimes_{h \A}$, and the fact that every 
nondegenerate Hilbert $\D$-module is a complemented submodule
of a direct sum of copies of the universal representation,
the result reduces to
proving that if (iii) holds for $K$, then it also 
holds for $K^{\gamma}$ for some cardinal
$\gamma$.   However this is easily seen from the
injectivity of
the Haagerup tensor product \cite{PS,BP1}, together
with the
operator space identification 
$K^{\gamma}_r \otimes_{h \A} V_k \cong
R_\gamma \otimes_h K_r \otimes_{h \A} V_k$, where 
$R_\gamma$ is the row Hilbert space of dimension $\gamma$.
That (i) is equivalent to (ii)
 follows easily from the universal properties
of $\D$-injectivity, and  \ref{upp}.  
For the `$\implies$' direction
take a completely contractive
$\A$-module map $T : V_1 \rightarrow \B$.
By \ref{up} we get a completely contractive
$\D$-module map $\D \otimes_{h\A} V_1 \rightarrow \B$.
Hence, by our hypothesis and $\D$-injectivity of $\B$,
there is a completely contractive
$\D$-module map extension $\tilde{T} : \D \otimes_{h\A} V_2
\rightarrow \B$.  Then $\tilde{T}$ restricted to
$V_2$ is a completely contractive
$\A$-module extension of $T$ to $V_2$. 
The other direction follows easily by showing that the
`closure of $\D \otimes_{\A} V_1$' in 
$\D \otimes_{h\A} V_2$ has the correct
universal property (in the remark after
\ref{up}).    \end{proof}
 
\vspace{4 mm}

\noindent {\bf Remarks.}
1) By symmetry, if we are concerned with right modules, the
analogous condition in (iii) would be in terms of spaces 
$V_k \otimes_{h \A} K$.  It is unnecessary to 
consider the dual space $K_r$.

\noindent 2) One may replace the `$\D$' by `$\A$' in 
condition (iii) and (iv) above, in the case that 
$\D = \C$.

\noindent 3) Let us say that a pair $(V_1,V_2)$ satisfying the 
equivalent conditions of the previous theorem, has the extension 
property.  For example,
if there is a completely contractive $\A$-module
projection $P : V_2 \rightarrow V_1$, then $(V_1,V_2)$ has 
this property for any such $\D$.
In particular if
$\D = \C$ and $V_2$ is a nondegenerate
Hilbert $\A$-module with submodule $V_1$, then 
$(V_1,V_2)$ has the extension property
if and only if the projection of $V_2$ onto $V_1$ is an 
$\A$-module map.  Thus  if $V_1$ is a fixed 
nondegenerate Hilbert $\A$-module, then $V_1$ is
orthogonally injective in the sense of \cite{MS}
if and only if $(V_1,V_2)$ has the extension property
whenever $V_2$ is a nondegenerate
 Hilbert $\A$-module containing $V_1$.
We remark that an $\A$-injective Hilbert module is
orthogonally injective, fairly clearly.
Clearly, \ref{buy}
is related to the topic of `commutant lifting'.

The following theorem may be viewed as a 
continuation of the
pretty Theorem 3.1 of  \cite{MS}; where  Muhly and Solel
give several Hilbert module characterizations of
 C$^*-$algebras.
Indeed the main ingredient of our proof below is the
equivalence of (i) and (v) below, which  
is part of their 
result.  We will therefore not prove this 
equivalence below.

We found that item (ii) was implied by (vi) or (vii),
so that it was natural to conjecture that it
alone characterized C$^*$-algebras.
After asking him this question, Christian Le Merdy 
kindly supplied a proof of it using Pisier's
$\delta$-norms \cite{LM}.  Later we found the proof below
using Muhly and Solel's result.  We will use
this fact in the next section.

\begin{theorem}
\label{cod}  The following are equivalent for an
operator algebra $\A$ with c.a.i.: \begin{itemize}
\item [(i)]  $\A$ is a C$^*-$algebra.
\item [(ii)]  For every completely contractive representation
$\pi : \A \rightarrow B(H)$, the commutant
$\pi(\A)'$ is selfadjoint.  
\item [(iii)]  $B(H)$ is (left) $\A$-injective for every 
Hilbert $\A$-module $H$.
\item [(iv)]  Every 
Hilbert $\A$-module $H$ is $\A$-injective.
\item [(v)]   For every nondegenerate
completely contractive representation
$\pi$ of $\A$ on a Hilbert
 space $H$, and every $\pi(\A)$-invariant 
closed subspace 
$K$ of $H$, $H \ominus K$ is $\pi(\A)$-invariant.
\item [(vi)]  $\C \otimes_{h\A} \C $ is completely isometrically
isomorphic to $\C$, as a $\C-\C$-operator bimodule.
\item [(vii)]   For every 
Hilbert $\A$-module $H$, the dilation 
$\C \otimes_{h\A} H$ is a Hilbert
space.
\item [(viii)] The canonical map from the $\C$-dilation
of $V_1$ to the $\C$-dilation
of $V_2$ 
is a complete isometry 
whenever $V_2$ is an $\A$-operator module
with closed submodule $V_1$.
  \end{itemize}
\end{theorem}

\begin{proof}    
By \ref{wit}, (i) implies (iii).
Clearly (iii) implies (iv), since $H$ is naturally a
complemented $\A$-submodule of $B(H)$.
That (iv) implies (v) is in \cite{MS}, since as we
said, an injective Hilbert module is orthogonally
injective, but in any case
the proof is immediate by extending 
the inclusion $i : K \subset H$ to a completely
contractive $\A$-module map $P \in B(H)$.  Clearly
$P$ is the projection onto $K$, and since it is an
$\A$-module map we obtain (v).  A similar idea shows that
(ii) implies (v); if $\pi$ is as in (v), and if $\theta$
is $\pi$ restricted to $K$, let 
$\rho = \theta \oplus \pi$, which is a 
representation of $\A$
on $K \oplus H$.   If $i$ is as above, then 
$$
T = \left[ \begin{array}{ccl}
0 & 0 \\
i & 0 
\end{array} \right]
$$   
commutes with $\rho$.  If (ii) holds, $T$ commutes 
with $\rho(\A)^*$, which easily gives (v). 
Thus (i)-(v) are all equivalent.

Clearly (i) implies (vi).
If we have (vi), 
and if $V$ is any $\C$-operator module, then
$\C \otimes_{h\A} V \cong \C \otimes_{h\A} \C \otimes_{h\C} V 
\cong \C \otimes_{h\C} V \cong V$ .
Taking $V$ to be a Hilbert $\C$-module shows 
(vii).  Assuming (vii), namely that
 $K = \C \otimes_{h\A} H$ is a Hilbert
space, write $i$ for the
canonical map $H \rightarrow K$ mentioned
in \ref{up}, and let $T : H \rightarrow H$ be the
identity map.  By \ref{upp}, there is a completely contractive
$\C$-module map $\tilde{T} : K \rightarrow H$ such
that $\tilde{T} \circ i = T$.  By Lemma \ref{su},
$i$ is a $\C-$module map.  Hence
$i$ is onto, which shows that 
$\C \otimes_{h\A} H \cong H$.
Hence, by the universal property
\ref{upp}, given an $\A$-submodule $H$ of any
Hilbert $\A$-module $K$ as in (v), the
inclusion map $i : H \rightarrow K$
is a $\C$-module map; and so we see that (v) holds.
Thus (i)-(vii) are all equivalent.

Clearly (i) implies (viii), whereas (viii) implies 
(iii) by \ref{buy}.
\end{proof}

\vspace{5 mm}

\noindent {\bf Remarks.}  If $H$ is a nondegenerate
Hilbert $\A$-module
the proof above shows that $\C \otimes_{h\A} H$ is a 
Hilbert space if and only if 
$\C \otimes_{h\A} H \cong H$.  As in `(vii) $\implies$ (v)' 
above, this implies that $H$ is an orthogonally 
injective module
in the sense of \cite{MS}, and also that the 
commutant in $B(H)$ of the associated representation of
$\A$ on $H$ is selfadjoint.  The converse is 
not true however.
Simple calculations show in
the case where
$\A$ is the disk algebra, then the only
Hilbert modules with
$\C \otimes_{h\A} H \cong H$, are one-dimensional. 
 In the case 
when $\A = {\mathcal T}(2)$, the upper triangular
$2 \times 2$ matrices, these modules coincide with 
the Hilbert $\ell^\infty_2$-modules - in other words the
nilpotent part of the action vanishes.

In \S 5 of
\cite{Pc}, Paulsen shows that if $H$ is a Hilbert module
over the disk algebra $\A = A({\mathbb D})$ associated 
with a coisometry then $B(H)$ is $\A$-injective.  Muhly
and Solel show in \cite{MS} 
that these $H$ are the `orthogonally 
injective Hilbert modules'.  This class of modules
coincides also with the 1-injective Hilbert 
$\A$-modules.

Notice that in (viii) one may replace $\C$ by any
C$^*-$algebra generated by $\A$.  

Finally, one can check that (ii) is equivalent to
 the universal representation of $\A$ having
 selfadjoint commutant.
 
\section{Morita equivalence of operator algebras}

In this section $\A$ and $\B$ are operator algebras with
c.a.i.   We refer the reader to \cite{BlFTM1,BMP} if 
further background for this section is needed.  For the
basic theory of Morita equivalence and strong 
Morita equivalence of C$^*-$algebras we refer the reader
to \cite{Ri1,Ri2,La}.

We begin with a brief discussion of `weak Morita
equivalence'.  
This is mostly independent of the rest of this section,
and the reader could skip to 5.2, if desired.   
Loosely speaking, this means that two operator 
algebras have `the same' Hilbert space
representations.   More precisely, we say that
$\A$ and $\B$ are weakly Morita equivalent
if the categories
$_{\A}HMOD$ and $_{\B}HMOD$ are naturally isometrically 
equivalent\footnote{That is, if there exist contractive
functors $F : $ $_{\A}HMOD \rightarrow $ $_{\B}HMOD$
and $G : $ $_{\B}HMOD \rightarrow $ $_{\A}HMOD$, such that 
$F G \cong Id$ and
$G F \cong Id$ naturally isometrically.}.
It is not hard to show that for C$^*-$algebras
$\C$ and $\D$, weak Morita
equivalence coincides with what was called `Morita 
equivalence'\footnote{In recent years
we have heard the term
`weak Morita equivalence' being used for Rieffel's 
`Morita equivalence of C$^*-$algebras' (as opposed to
his `strong Morita equivalence').}
in \cite{Ri1}
We note that it is folklore that the 
latter happens if and only if
there is a Hilbert space $H$ such that
$e(\C) \bar{\otimes} B(H) \cong e(\D) 
\bar{\otimes} B(H)$ *-isomorphically,
where $e(\C)$ is the enveloping 
von Neumann algebra of $\C$.

Henceforth we reserve the symbols $\C$ and $\D$ for the
maximal C$^*-$algebras generated by $\A$ and $\B$ respectively.

\begin{proposition}
 \label{wme}  If $\A$ and $\B$ are weakly Morita 
equivalent operator algebras then:
\begin{itemize}
\item [(i)]  If $\A$ is a C$^*-$algebra then so is $\B$.
\item [(ii)]  $\C$ is weakly Morita equivalent to $\D$.
\end{itemize}
\end{proposition}
 
\begin{proof}  Suppose that $F : \; _{\A}HMOD
\rightarrow \; _{\B}HMOD$ is an equivalence functor.
For $H \in $ $_{\A}HMOD$, we have
$_{\C}B(H)$ is a subalgebra of $_{\A}B(H)$.  By an
obvious argument (see for example Lemma
2.2 in \cite{BlFTM1}),
 the map $T \mapsto F(T)$ from $_{\A}B(H)
$ to $_{\B}B(H)$ is a
isometric homomorphism.  Hence its
restriction to the C$^*-$algebra $_{\C}B(H)$ is a
*-homomorphism, and consequently maps into
$_{\D}B(F(H))$.

From this we see that if $\A = \C$, then $_{\B}B(H)$ is
a C$^*-$algebra for all Hilbert $\B$-modules.   By the
implication `(ii) $\implies$ (i)' in  
Theorem \ref{cod}, we see that $\B$ is a C$^*-$algebra. 
 
Now suppose that  $H_1, H_2 \in $ $_{\A}HMOD$ and $T \in
\; _{\C}B(H_1,H_2)$.
Let  $H = H_1 \oplus H_2$, and let
 $i_k$ and $q_k$ be, respectively,
 the inclusions and projections between the $H_k$ and $H$.
Thus $q_k \circ i_k = Id_{H_k}$, so that
$F(q_k) F(i_k) = Id_{F(H_k)}$.  From \ref{su} it 
follows that $i_k, q_k, F(q_k)$
 and $ F(i_k)$ are $\C-$module maps.
By the first part, $F(i_2 T q_1)$ is a $\C-$module map.
Thus $F(T) = F(q_2) F(i_2 T q_1) F(i_1)$ is a
$\C-$module map.

We have shown that  $F$ restricts to a functor from
$_{\C}HMOD$ to $_{\D}HMOD$.  Similarly for $G$, and
now the
category equivalence is clear.  Note that the
natural transformation maps are unitary and commute
with the action of the operator algebra, and hence
also commute with the action of the generated C$^*-$algebra.
Thus we have (ii).  
\end{proof}

\vspace{5 mm}

We refer the reader to \cite{BMP} for the 
definition of {\em strong Morita equivalence} of
operator algebras $\A$ and $\B$.  Loosely,
it is defined in terms of two
operator bimodules $X$ and $Y$, which possess certain
pairings $(\cdot) : X \times Y \rightarrow \A$ and
$[\cdot] : Y \times X \rightarrow \B$.  
The tuple $(\A,\B,X,Y,(\cdot),[\cdot])$ is called a 
{\em strong Morita context} (see \cite{BMP} 
Definition 3.1).  Here we shall usually simply 
write $(\A,\B,X,Y)$.  This generalizes 
C$^*-$algebraic strong Morita equivalence \cite{Ri2}.
If $\A$ and $\B$ are C$^*-$algebras it turns out that 
$X$ may be taken to be the conjugate
bimodule of $Y$ (or equivalently, $X = Y^*$ in the linking
C$^*-$algebra \cite{La}).

\begin{definition}
\label{subco}
\end{definition}
\begin{itemize}
\item [(i)] Suppose 
that $\E$ and $\F$ are strongly Morita equivalent 
C$^*-$algebras, and that $Z$ is an $\F-\E$-strong
Morita equivalence bimodule, and that $W = \bar{Z}$ is the
conjugate $\E-\F-$bimodule of $Z$.  Then we say that
$(\E,\F,W,Z)$ is a {\em C$^*-$Morita context},
or {\em C$^*-$context} for short.
\item [(ii)] Suppose that $\A$ and $\B$ are operator algebras
with c.a.i., and suppose that $\E$ and $\F$ are C$^*-$algebras
generated by $\A$ and $\B$ respectively.  Suppose that 
$(\E,\F,W,Z)$ is a C$^*-$Morita context,  
$X$ is a closed $\A-\B$-submodule of $W$, and
that $Y$ is a
closed $\B-\A$-submodule of $Z$.  Suppose further
that the natural 
pairings  $Z \otimes W \rightarrow \F$ and 
$W \otimes Z \rightarrow \E$ restrict to
maps $Y \otimes X \rightarrow \B$, and $X \otimes Y
\rightarrow \A$, both with dense range.  Then we say that
$(\A,\B,X,Y)$ is a {\em subcontext} of $(\E,\F,W,Z)$.    
If, further,
 $\E$ and $\F$ are the maximal C$^*-$algebras of 
$\A$ and $\B$ respectively, then we shall
say that $(\A,\B,X,Y)$ is a {\em maximal
subcontext.}  Similarly, a {\em minimal
subcontext} occurs when $\E$ and $\F$ are the
C$^*-$envelopes of $\A$ and $\B$.
 \item [(iii)]  A subcontext 
$(\A,\B,X,Y)$ of a C$^*-$Morita context
$(\E,\F,W,Z)$ is said to be  {\em left dilatable}
if $W$ is the left $\E$-dilation of $X$, and 
$Z$ is the left $\F$-dilation of $Y$.  In this case
we say that $\A$ and $\B$ are {\em left strongly 
subequivalent}.  We also say that $X$ and $Y$ are
(left) subequivalence bimodules, and that 
$(\A,\B,X,Y)$ is a {\em left subequivalence context}.
\end{itemize}

There is a similar definition and symmetric theory
where we replace
the words `left', by  `right' or
`two-sided'.  Generally, we shall omit the word
`two-sided' and simply refer, for example to
`strong subequivalence'.  

In order to come to grips with these definitions, we  
proceed with several observations and examples: 

\noindent {\bf Remarks.}  Note that (ii) implies that
$X$ and $Y$ are nondegenerate operator bimodules over
$\A$ and $\B$. This is because $W$ and $Z$ are automatically
nondegenerate (see 1.5 in \cite{La}), and any c.a.i for 
an operator algebra is also a c.a.i. for any C$^*-$algebra it
generates.

Write $\Li$ for the 
set of $2 \times 2$ matrices 
$$
\left[ \begin{array}{ccl}
a & x \\
y & b
\end{array} \right] 
$$
with $a \in \A, b \in \B, x \in X, y \in Y$.  Write
 $\Li'$ for the same set, but with entries from the
C$^*-$context $(\E,\F,W,Z)$.  It is well known 
(see \cite{La}) that
$\Li'$ is canonically
a C$^*-$algebra, called the `linking
C$^*-$algebra' of $Z$, or of $(\E,\F,W,Z)$.  
Saying that 
$(\A,\B,X,Y)$ is a subcontext of $(\E,\F,W,Z)$
is almost equivalent to saying that 
$\Li$ is a closed subalgebra of $\Li'$.  We say `almost',
because the latter condition does not imply 
the statement in (ii) about `dense range'.  In any 
case it is clear that  a
subcontext gives a {\em linking operator algebra} $\Li$.
Clearly $\Li$ has a c.a.i.  We shall see
that $\Li$ generates $\Li'$ as a C$^*-$algebra.

If $(\A,\B,X,Y)$ is a subcontext of
$(\E,\F,W,Z)$, and if $\A$ and $\B$ are unital, then
the pairings in (ii) having dense range is equivalent to
(as in Proposition 3.3 of \cite{BMP})
these pairings being  onto, and hence
$(\A,\B,X,Y)$ is a `c.b.-Morita
context' in the sense of \cite{BMP} Definition 3.1.
However, we are mainly interested in
when a subcontext is
a strong Morita context.  
 
Note that in (iii) we are in the situation where the
canonical map from the operator module
($X$ or $Y$)  into its dilation ($W$ or $Z$)
is a complete isometry.
We shall see later that there are some simple tests 
for when a subcontext is left dilatable.  

Finally, we remark that we 
showed in \cite{BMN} that strong Morita
equivalence implies (two-sided) strong 
subequivalence, and moreover the 
implicit subcontext may be taken to be minimal 
(or maximal).  Together with
K. Jarosz we have found a simple example \cite{BJ}
of a closed subalgebra $\A$ of the disk algebra, 
giving a strong subequivalence (which is a
minimal sucontext, and is of the type discussed 
in Example 3 below) which is
not a strong Morita equivalence.  Thus strong 
subequivalence is genuinely a new notion.  We 
shall see however, that strong Morita equivalence is
the same as strong subequivalence 
via a {\em maximal} subcontext.
Hopefully the distinctions will be illuminated
more clearly as we go along.  We also refer to 
\cite{BJ} for further, and very concrete,
 illumination of these notions.

\vspace{4 mm}

\noindent {\bf Examples 1.)} The `dense range' 
condition in (ii)
is not implied by the dilation condition
in (iii).  Indeed if $\A = {\mathcal T}(2)$,
$\B = {\mathbb C}$, $Y = R_2$, and $X = {\mathbb C}
[1 \; \; 0]^t$, and if $\E = M_2$ (the C$^*-$envelope of 
$\A$) and $\F = \B$, then it's easily seen that
(ii) and (iii) hold with the exception of the
pairing $X \otimes Y \rightarrow \A$ having dense range.
This example is interesting in that in this case
the C$^*-$envelopes of $\A$ and $\B$ are strongly
Morita equivalent with equivalence bimodules being the 
minimal C$^*-$dilations of $X$ and $Y$ above, but the 
maximal C$^*-$algebras of $\A$ and
$\B$ are not Morita equivalent in any sense. 

{\bf 2.)}  Another interesting example of 
subcontexts comes from
example 8.2 in \cite{BMP} (see also
6.9 in \cite{Bhmo}, where $\A = A({\mathbb D})$
 is the disk algebra,
and we find $\A$-operator modules $X, Y$ 
such that $(\A,\A,X,Y)$ and $(\A,\A,\A,\A)$ are two 
different (two-sided) dilatable subcontexts
of $(\E,\E,\E,E)$, where $\E = C({\mathbb D})$.  
Hence one cannot hope in general to recover $X,Y$ from 
the data of $\A,\B$ and the containing 
C$^*-$context $(\E,\F,W,Z)$.  In example 8.3 of 
\cite{BMP} we discussed another subcontext coming from 
matrix algebras of analytic functions, which
is not dilatable.  

{\bf 3.)}  Related to the Example 2, let $\A$ be
an operator algebra with identity of norm 1, 
let $\E$ be a 
C$^*-$algebra generated by $\A$, and choose $x \in
\E \setminus \A$, with $x$ invertible in $\E$, such that
$x^{-1}\A x$ generates $\E$ (such as is the case when
$x \in \A'$).
Then $(\A,x^{-1}\A x,\A x , x^{-1} \A)$ is a
subcontext of the `identity context'
$(\E,\E,\E,\E)$.  One may quite
easily write down conditions on $x$ ensuring that
this subcontext is left dilatable.   
For example, suppose
that $\Omega$ is a compact Hausdorff space, and $\A$ is a
uniform algebra on $\Omega$ 
(containing constants and separating
points).  Let $P = \{ |g| : g \in \A \} 
\subset C(\Omega)_+$.  
Choose a strictly positive  function $f$ on $\Omega$,
 such that $f, f^{-1} \in P$, or equivalently: 
$f \in P \cap P^{-1}$.  Then it is easy to see from the
Stone-Weierstrass theorem that
$(\A,\A , \A f , f^{-1} \A)$ is a (two-sided) dilatable 
subcontext
of the `identity context' of $C(\Omega)$.  In fact this 
is true
under much less restrictive conditions on $f$.
It appears to be an interesting 
function algebra question to characterize when such subcontexts
are strong Morita contexts 
(see \cite{BJ} for more details).  
For example, it is easy to see
that they always are, if 
$f \in Q$ or $f \in Q^{\bar{}}$ (the uniform closure of $Q$), 
where $Q = 
\{ |k| : k , k^{-1} \in \A \} \subset P \subset C(\Omega)$.
Note that $P \cap P^{-1} = Q$  if  $\A$
is the disk algebra, say, and $\Omega = \Di$.
This is because if $f = \vert g \vert,
f^{-1} = \vert h \vert, g, h \in \A$, then 
$\vert gh \vert = 1 $, and hence
by the maximum modulus theorem
$gh$ is constant, and hence $g$ is invertible in $\A$.

 \vspace{5 mm}

\begin{proposition}   
\label{subcont}
If 
$(\A,\B,X,Y)$ is a subcontext of a C$^*-$Morita context
$(\E,\F,W,Z)$,
then 
\begin{itemize}
\item [(i)]  $X$ and $Y$ generate
$W$ and $Z$ respectively as left operator
modules.  
 So, for example, $W$ is the smallest closed 
left $\E$-submodule of $W$ containing $X$.  
Similar assertions hold as right operator
modules, by symmetry.  
\item [(ii)] The 
linking operator algebra $\Li$ generates the 
linking C$^*-$algebra $\Li'$ of $(\E,\F,W,Z)$.     
\item [(iii)]   
If $\A$ {\em or}  $\B$ is a C$^*-$algebra, then
$(\A,\B,X,Y) = (\E,\F,W,Z)$.
\end{itemize}
\end{proposition}
  
\begin{proof}  It is easy to see that (ii) and (iii)
follow from (i).
We shall simply show that 
$X$ generates
$W$ as a left $\E$-operator module.
Since the pairing $
[\cdot] : Y \otimes X \rightarrow \B$
has dense range, we can pick a c.a.i. for $\B$ which is a
sum of terms of the form $[y,x]$, for $y \in Y, x \in X$.
This c.a.i. is also one for $\F$, and hence
sums of
terms of the form $w [y,x]$, for $y \in Y, x \in X,
w \in W$ are dense in $W$.  However,
$w [y,x] = (w,y) x \in \E X$ (where $(\cdot)$ is the
other pairing).  So $X$ generates
$W$ as a left $\E$-operator module.  
\end{proof}

\begin{theorem}
\label{sep}  If $(\A,\B,X,Y)$ is a strong Morita
context 
which is a subcontext of a C$^*-$Morita context
$(\E,\F,W,Z)$,
 then it is a dilatable subcontext. 
\end{theorem}

\begin{proof}  By the previous result,
$X$ and $Y$ generate
$W$ and $Z$ respectively as left operator
modules.
Thus we have a
complete contraction $\E \otimes_{h\A} X
\rightarrow W$ with dense range.  On the other hand
$$
W \cong W \otimes_{h\B} \B \cong
W \otimes_{h\B} Y \otimes_{h\A} X
\cong (W\otimes_{h\B} Y) \otimes_{h\A} X \; \; .
$$
However, the pairing $(\cdot)$ determines a
complete contraction $W \otimes_{h\B} Y
\rightarrow \E$, and so we obtain a complete
contraction $W \rightarrow \E \otimes_{h\A} X$.
One easily checks that the composition of these maps
$$\E \otimes_{h\A} X
\rightarrow W  \rightarrow \E \otimes_{h\A} X$$
is the identity, from which it follows
 that $W \cong \E \otimes_{h\A} X$.
Similarly $Z$ is the dilation of $Y$.
\end{proof}

\begin{theorem}
\label{suf}  If $(\A,\B,X,Y)$ is a
left dilatable
maximal subcontext of a C$^*-$context, then
$\A$ and $\B$ are strongly Morita equivalent
operator algebras,
and $Y$ is a strong $\A-\B$-Morita
equivalence bimodule, with dual module $X$.  Indeed,
it also follows that $(\A,\B,X,Y)$ is a (strong)
 Morita context.
Conversely, every strong Morita equivalence
of operator algebras occurs in this way.  
That is, every strong Morita context is a
left dilatable
maximal subcontext of a C$^*-$Morita context.
\end{theorem}

\begin{proof} 
If $\C$ and $\D$ are as usual the maximal 
C$^*-$algebras of $\A$ and $\B$ respectively, and if
$(\A,\B,X,Y)$ is a left dilatable
subcontext of $(\C,\D,W,Z)$ then, using Lemmas
\ref{sub} and \ref{up}, we
have $$
Y \otimes_{h\A} X
\subset \D \otimes_{h\B} (Y \otimes_{h\A} X)
\cong Z \otimes_{h\A} X
\cong (Z \otimes_{h\C} \C) \otimes_{h\A} X
\cong Z \otimes_{h\C} W \cong \D \; \;, 
$$
completely isometrically.
On the other hand, we have the canonical complete
contraction
$$
Y \otimes_{h\A} X \rightarrow \B \subset \D \; \; , $$
coming from the restricted pairings in (ii).
It is easy to check that the composition of the maps in these
two sequences agree.  Hence the canonical map
$Y \otimes_{h\A} X \rightarrow \B$ is a completely isometric
isomorphism.   Similarly, $X \otimes_{h\B} Y \cong \A$
completely isometrically.  Thus by the remark 
before Definition 3.6 in \cite{BMP}
(or see Definition 1.2 in \cite{BlFTM2} and the `sketch'
beneath it), $\A$ and $\B$ are
strongly Morita equivalent operator algebras.
 
The last statement
is proved in \cite{BMN}.
\end{proof}

\vspace{5 mm}

The last theorem may be viewed as a new characterization of
strong Morita equivalence of operator algebras.

\vspace{5 mm}

We next show that `strong subequivalence'
seems to have many of the nice implications of
 strong Morita equivalence (see Theorem 4.1 in \cite{BMP}
and the end of Chapter 3 there).  We intend to pursue in 
the near future exactly which other of the consequences
of strong Morita equivalence still carry over for this
weaker notion.   There is presumably also a theory of 
`sub-rigged' modules paralleling notions from 
\cite{Bhmo}, although we expect to lose some of
the rich features of rigged modules.

\begin{theorem}
\label{mo1}  Suppose that $(\A,\B,X,Y)$ is a
left dilatable
subcontext of a C$^*-$context $(\E,\F,W,Z)$.
Then $Y \cong $$_{\A}CB^{ess}(X,\A)$ and 
$ X \cong $$_{\B}CB^{ess}(Y,\A)$ completely
isometrically and as operator bimodules,
and $\A \cong $$_{\B}CB^{ess}(Y,Y)$ and 
$\B \cong $$_{\A}CB^{ess}(X,X)$ 
completely
isometrically and as operator algebras.
Moreover, the categories of ${\A}$-submodules
of $\E$-operator modules, and $\B$-submodules
of $\F$-operator modules, are (completely 
isometrically) equivalent.  This equivalence 
restricts to an 
equivalence of the categories of ${\A}$-submodules
of Hilbert $\E$-modules, and $\B$-submodules
of Hilbert $\F$-modules.   
\end{theorem}

\begin{proof}  Write $(\cdot,\cdot)$ and
$[\cdot,\cdot]$ for the pairings discussed in (ii)
of Definition \ref{subco}.    
Notice firstly, that there is a natural map
$Y \rightarrow \; _{\A}CB^{ess}(X,\A)$ coming 
from these
pairings.  Hence we get a sequence
$$
Y \rightarrow \; _{\A}CB^{ess}(X,\A) \; \subset \;
_{\A}CB^{ess}(X,\E) \;
 \cong \; _{\E}CB^{ess}(W,\E) \; \cong \; Z \; \; ,
 $$
where the second last map comes from \ref{bim}.
However, the composition of maps in this sequence agrees with
the inclusion of $Y$ in $Z$.  Hence the map
$Y \rightarrow \; _{\A}CB^{ess}(X,\A)$ is 
a complete isometry.
That this map is onto follows by the argument
of \cite{BMP} Theorem 4.1.
A similar proof shows that
$X \cong \; _{\B}CB^{ess}(Y,\B)$ as operator bimodules,
and that $\B \cong
\; _{\A}CB^{ess}(X,X)$ and $\A \cong \; _{\B}CB^{ess}(Y,Y)$
(completely isometrically) as operator algebras.
Define $F(V) = \; _{\A}CB^{ess}(X,V)$ and $G(U) =
\; _{\B}CB^{ess}(Y,U)$, we will show that  $F$ and $G$ are
are completely contractive
equivalence functors between the category 
operator $\A$-submodules of $\E$-operator modules,
and the category 
of $\B$-submodules of $\F$-operator
modules, which compose (up to natural completely
isometric isomorphism) to the
identity functor. 

If $V$ is an $\E$-operator module, then by \ref{bim} and
\cite{Bna} Theorem 3.10, we have
$$
F(V) = \; _{\A}CB^{ess}(X,V) \cong \;
_{\E}CB^{ess}(W,V) \cong
Z \otimes_{h\E} V \; \; .
$$
Moreover,
this, together with the corresponding result for $G$,
shows that $G(F(V)) \cong W \otimes_{h\F} Z \otimes_{h\E} V
\cong \E \otimes_{h\E} V
 \cong V$.

For a general $\A$-operator module $V$,
there is a canonical complete
contraction $\rho_V : V \rightarrow
G(F(V)) = \; _{\B}CB^{ess}(Y,_{\A}CB^{ess}(X,V))$ 
given by
$(\rho_V(v)(y))(x) = (x,y) v$, for
$v \in V , y \in Y , x \in X$.
Suppose that  $V$ is an $\A$-operator module, 
and that $V'$ is a
$\E$-operator module containing $V$.  Then we get
the following sequence of complete contractions
$$V \rightarrow G(F(V)) = \; 
_{\B}CB^{ess}(Y,_{\A}CB^{ess}(X,V)) \subset
\; _{\B}CB^{ess}(Y,_{\A}CB^{ess}(X,V')) \cong V' \; \; .
 $$
The first map here is
$\rho_V$.  The composition of maps in this sequence
is the inclusion map, and so $\rho_V$ is a complete
isometry.  To show that $\rho_V$ is onto in the unital 
case is a simple
exercise in algebra.  In the nonunital case, to
show that $\rho_V$ is onto, one may use an
 argument similar
to those in the proof of \cite{BMP} Theorem 4.1 showing that
the maps there are onto.  That Hilbert modules in these 
categories are taken
by this equivalence to Hilbert modules follows
easily from the observations above.
\end{proof}

\vspace{5 mm}

We recall (see \cite{MS} for example) that a 
Shilov Hilbert module is an $\A$-submodule of a 
Hilbert module over the C$^*-$envelope of $\A$.
As a consequence it follows that 
minimal subequivalence of two operator algebras
implies a (weak) equivalence 
between the subcategories of Shilov Hilbert modules.
We are led to
propose the following definition:

\begin{definition}  
\label{subM}  We say that operator algebras 
$\A$ and $\B$ are  (two-sided) 
{\em minimally subequivalent} if they are (two-sided) 
strongly subequivalent, and 
the C$^*-$algebras in the containing C$^*-$context 
are the C$^*$-envelopes of $\A$ and $\B$.
\end{definition}

A similar definition pertains where we replace
the word `two-sided' by  `left' or   `right'.  Notice that 
there is no need to
define `maximally subequivalent', since this would
coincide with strong Morita equivalence, by Theorem
\ref{suf}.   Strong Morita equivalence implies
minimal subequivalence by \cite{BMN}.   However, we have  
examples to 
show that the converse is false:
indeed two-sided minimal subequivalence is a  
weaker notion than strong Morita equivalence.

We now show how `strong subequivalence' can arise,
by discussing some equivalent conditions for a 
subcontext $(\A,\B,X,Y)$ of $(\E,\F,W,Z)$
to be left dilatable, or equivalently, for
$W \cong \E \otimes_{h\A} X$ and $Z \cong
\F \otimes_{h\B} Y$.
As we saw in 
\ref{subcont}, the definition of a 
subcontext already implies 
that $[\E X]^{\bar{}} = W$ and $[\F Y]^{\bar{}} = Z$.

Theorem \ref{eno} or \ref{bim} tells us
that $W \cong \E \otimes_{h\A} X$ 
 is equivalent to  
the fact that
$_{\E}CB^{ess}(W,H) \cong \; _{\A}CB^{ess}(X,H)$ completely 
isometrically for all
Hilbert $\E$-modules $H$.  Indeed the Hilbert
space of the universal representation would suffice.
From C$^*-$module theory 
(see \cite{Bna} for background) we have that
$_{\E}CB^{ess}(W,H) \cong Z \otimes_{h\E} H$.
The last space is a Hilbert column space,
whose norm we can completely describe:
namely $\Vert \sum_k z_k \otimes \zeta_k \Vert^2 
= \sum_{k,j} \langle \langle z_k \vert z_j \rangle \zeta_j ,
\zeta_k \rangle \; .$   Here the inside $\langle \cdot \rangle$
is the $\E$-valued inner product on $Z$.
A similar formula gives the
matrix norms (see \cite{BMP} Lemma 2.13). 
The restriction map $_{\E}CB^{ess}(W,H) 
\rightarrow \; _{\A}CB^{ess}(X,H)$
may thus be rewritten as the map $R : Z \otimes_{h\E} H
\rightarrow \; _{\A}CB^{ess}(X,H)$, 
given by $R(z \otimes \zeta)(x) = (x,z) \zeta$.
By \ref{eno}, we
need to check that $R$ is a complete isometry
($R$ is onto by similar considerations to those in the 
proof of the
previous theorem).  
If $_{\A}CB^{ess}(X,H)$ is known to be a Hilbert column
space (which is the case, say, if we know that $X$ is
a left $\A$-rigged module \cite{Bhmo}), 
then we need only check that $R$ is an isometry, or 
equivalently that 
$$ \sum_{k,j} \langle \langle z_k \vert z_j \rangle \zeta_j ,
\zeta_k \rangle \; \leq \; \sup \{ \Vert [ \sum_k (x_{ij},z_k) \zeta_k ] 
\Vert^2 \} \; \; ,
$$
whenever $z_1 , \cdots , z_n \in Z$, $\zeta_1 , \cdots , 
\zeta_n \in H$,
where the supremum is taken over all sized matrices $[x_{ij}]$ of 
norm $1$ with entries in $X$.  

The second part of \ref{eno} gives 
another condition which is equivalent to the above, and
which may be easier to check in a concrete example: 
namely that the canonical completely contractive map
$S : H_r \otimes_{h\A} X \rightarrow 
H_r \otimes_{h\E} W$, is an isometry.  In this case
$H_r \otimes_{h\E} W$ is a row Hilbert space, so that 
if $S$ is an isometry then it is automatically a 
`row-isometry'.  

By symmetry, the subcontext $(\A,\B,X,Y)$ of
$(\E,\F,W,Z)$ is right dilatable if and only if
$X \otimes_{h\B} K \cong W \otimes_{h\F} K$ 
and $Y \otimes_{h\A} H \cong Z \otimes_{h\E} H$
isometrically, via the canonical maps, for all 
Hilbert $\E$-modules $H$ and all Hilbert 
$\F$-modules $K$.  Clearly, the condition 
$Y \otimes_{h\A} H \cong Z \otimes_{h\E} H$
for example, is equivalent to saying that the $\otimes_{h\A}$-norm
on $Y \otimes H$ equals:
$$
\Vert \sum_{j=1}^n y_j \otimes \zeta_j \Vert^2
= \sum_{k,j} \langle \langle y_k \vert y_j \rangle \zeta_j ,
\zeta_k \rangle \; ,   \eqno{(\dagger)}
$$
for $y_1 , \cdots , y_n \in Y$, $\zeta_1 , 
\cdots , \zeta_n \in H$.
Thus `right dilatability' is equivalent to 
saying that the induced functors
$F_Y = Y \otimes_{h\A} -$ and $G_X = X \otimes_{h\B} -$
coincide, on the categories of Hilbert $\E-$ and 
$\F-$modules, with the weak Morita equivalence induced by
$Z$ and $W$ of these
categories.  Summarizing:

\begin{corollary}
A subcontext $(\A,\B,X,Y)$ of a strong Morita 
equivalence of $\E$ and $\F$, 
is right dilatable
if and only if the induced functors $F_Y$ and $G_X$
give back the original weak Morita equivalence
as explained above.
This is equivalent to ($\dagger$) holding
for all Hilbert
$\E$-modules $H$,
and the analogous formula for 
$X \otimes_{h\B} K$ holding for all 
 Hilbert $\F$-modules $K$.
This is also equivalent to the canonical maps
$X \otimes_{h\B} K \rightarrow W \otimes_{h\F} K$
and $Y \otimes_{h\A} H \rightarrow  Z \otimes_{h\E} H$
being row isometric, where $H$ and $K$ are 
the universal representations of $\E$ and $\F$
respectively.
\end{corollary}

\begin{proof}
Only the last statement still needs a word of proof,
  and this is similar to the proof that (iv) implies
(iii) in \ref{buy}.
\end{proof}

A simple modification of the
 first of our examples of subcontexts shows
again that the dense range condition in 
the definition of a subcontext is necessary for
the corollary to hold.  Without it
one may have $F_Y$ and $F_X$ giving the same
weak Morita equivalence between $_{\E}HMOD$ and
$_{\F}HMOD$ as $F_Z$ and $F_W$, without
$\A$ and $\B$ being 
weakly Morita equivalent.

This ends our discussion of subcontexts.
A natural question is if there is a comparable
theory of quotient Morita contexts: 
We end this section with a simple but important
observation, which for some reason we overlooked when
writing \cite{BMP}.
Suppose that $\A$ and $\B$ are
strongly Morita equivalent operator algebras,
and that $(\A,\B,X,Y)$ is the associated Morita context.
Suppose that $\pi : \A \rightarrow \D$ is a completely
contractive homomorphism into an operator algebra
$\D$.  Then if $\E$ is the closure of the range
of $\pi$, then there exists a natural
Morita context $(\E,\F,P,Q)$, which we shall call the
{\em pushout} of $(\A,\B,X,Y)$ along $\pi$, which one may
construct as follows.  Suppose that $\E$ is a nondegenerate
subalgebra of $B(H)$.  Then $\pi$ may be viewed as a
representation of $\A$ on $H$.  The original Morita
context gives rise to
 a Hilbert space $K = Y \otimes_{h \A} H$, as in \cite{BMP}
Theorem 3.10, and an induced representation $\theta$ of $\B$
on $K$.  Indeed, since $\A$ is strongly Morita equivalent
to its linking algebra, we obtain an induced
completely contractive
representation $\rho$ of
the linking algebra $\Li$ of the Morita context
 $(\A,\B,X,Y)$ (see \S 5 of
\cite{BMP}) on a Hilbert space
$N = S \otimes_{h \A} H$, where $S$ is the bimodule
for the equivalence of $\A$ and its linking algebra.
In fact $S = \A \oplus_c Y$, with notation as in
\cite{Bhmo}.  By the associativity relations
on p. 411 of that paper, we see that
$N = H \oplus K$ .  Indeed we have recaptured the
`obvious' representation of $\Li$ on
$H \oplus K$, namely the one which is determined by
$$
\rho \left(  \left[ \begin{array}{ccl}
a & x \\
y & b
\end{array} \right] \right) \left[ \begin{array}{cl}
\zeta \\
y' \otimes \xi
\end{array} \right]
\; \; = \; \;
\left[ \begin{array}{cl}
\pi(a) \zeta + \pi((x,y')) \xi \\
y \otimes \zeta + by' \otimes \xi
\end{array} \right]
\; \;.
$$
Here $\zeta, \xi \in H$.
The image under $\rho$ inside
$B(H \oplus K)$, of the
four corners of the linking algebra $\Li$, gives
a Morita context implementing a strong Morita equivalence
of $\E$ and the operator algebra which is the
closure of $\theta(\B)$.  This is because the complete
quotient conditions in the definition of strong
Morita equivalence (3.1 of \cite{BMP}), may be
checked by the lifting
criterion of 2.11 of \cite{BMP}.  This criterion,
loosely speaking,
is in terms of writing the c.a.i. of the algebras
in terms of elementary tensors of Haagerup
norm $< 1$.  However, if  $\{e_\alpha\}$
is a c.a.i. for $\A$, then $\pi(e_\alpha)$ is a
c.a.i. for $\E$; and the associated aforementioned
elementary tensors in $X \otimes Y$,
are taken, via the completely contractive
$\rho$,
 to elementary tensors in the
new context, of Haagerup
norm $< 1$.  Similarly for a c.a.i. for $\B$.

A special case of the pushout occurs when $\pi$ is the
quotient homomorphism associated with a closed 2-sided
ideal.  This case was been studied recently,
independently, and in much greater detail, in
 \cite{MS2}.
For C$^*-$algebras of course this is not
a `special case', but an equivalent formulation,
and was worked out in \cite{Ri3}.  Note that a 
pushout or `quotient Morita context' of a C$^*-$context
is again a C$^*-$context, because completely
contractive homomorphisms on C$^*-$algebras are
*-homomorphisms.  On the other hand, it 
might be interesting to determine when the pushout
of a strong Morita equivalence of operator algebras
is a C$^*-$context.

\section{The linking algebras of an
operator space or operator bimodule.}

We turn to another interesting connection between
the maximal C$^*-$algebra of an operator algebra,
and Morita equivalence/induced representations.  It is also 
interesting in that it gives rise to a class of
examples of Hilbert modules and C$^*$-modules which may
be associated to any operator space or operator module.  
It may also be viewed as a 
generalization of Example \ref{t2}, along an avenue opened
up by C. Zhang in \cite{Zh}.  He however was studying 
different questions, and was interested
in the C$^*-$envelope.
For clarity we will give the idea first in  the operator space
case, and then later discuss the more general operator 
bimodule case.

 Let $X$ be any operator space.  Assume that
$X \subset B(K,H)$ completely isometrically.  We form an
operator system ${\mathcal S}$ consisting of matrices
$$
\left[ \begin{array}{ccl}
\lambda_1 I_H  & x \\
y^* & \lambda_2 I_K
\end{array}
\right] \; \;
$$
where $x,y \in X , \lambda_1 , \lambda_2 \in {\mathbb C}$.
In this section $X^*$ will always mean the space
of adjoints, not the dual space.
A simple modification of Lemma 7.1 of
\cite{P}, or Theorem \ref{pr} below,
 shows that ${\mathcal S}$ is independent of the 
particular
$H,K$ chosen, up to completely isometric isomorphism.
Setting the 2-1 corner equal to $0$, gives a
unital operator algebra $\U(X)$, 
which only depends on the operator
space structure of $X$.  Let 
$\U_d(X)$ be the subalgebra
with repetition on the diagonal, and let $\A = \U(X), \A_d =
\U_d(X)$.  Let $\Li(X) = C^*_{max}(\A)$
and $\Li_d(X) = C^*_{max}(\A_d)$.   Given a
completely
contractive unital representation
$\pi$ of $\A_d$ on a Hilbert space $N$, the 
restriction of
$\pi$ to the $1-2$ corner gives a completely
contractive linear map $
\phi : X \rightarrow B(N)$.  Since $[\phi(X)N]^{\bar{}}$
and $(\cap_{x \in X} \ker \phi(x))^\perp$ are nontrivial 
complementary subspaces of $N$, we obtain a
a nontrivial
decomposition $N = H \oplus K$
say, with respect to which $\phi$ may be viewed as
a map $X \rightarrow B(K,H)$.  
Using the
aforementioned
modification of the
result in \cite{P}, $\phi$ may be `extended' to
a contractive unital representation
$\tilde{\pi}$ of $\A$ on $N$, which is also an
extension of $\pi$.  It follows from this that
$\Li_d(X)$ is a unital C$^*-$subalgebra of $\Li(X)$.

Thus there are  1-1 correspondences
between the following classes:
1) completely
contractive linear maps $X \rightarrow B(K,H)$,
 2) unital completely
contractive representations of $\U(X)$ on a
Hilbert space $N (= H \oplus K)$,
and  3) unital *-representations of $\Li(X)$
on $N (= H \oplus K)$; and moreover
one may use 
$\U_d(X)$ and $\Li_d(X)$
instead of $\U(X)$ and $\Li(X)$ in 2) and 3).
If $X$ is a maximal operator
space, then one may remove the words `completely' in
1) and 2) above.

The canonical projections $e_1, e_2$
in $\A$ give a decomposition
of $\Li(X)$ and $\Li_d(X)$ as $2 \times 2$ matrices.
By computations similar to the
exercise \ref{t2}, which are done more explicitly in
\cite{Zh},
one sees that $e_1 \Li(X) e_1$ is the closed linear span
of $e_1$ and terms
of the form $(x y^*)^n$, where $n \in {\mathbb N}$ and
$x,y \in X$, and the products here are
 with respect to $\Li(X)$.  Write
$\C$ 
or C$^*_{max}(XX^*)$ for $e_1 \Li(X) e_1$,
write $\D$ or C$^*_{max}(X^*X)$ for $e_2 \Li(X) e_2$, and 
write $W$ for the
`1-2 corner' $e_1 \Li(X) e_2$.  We call $W$ the
{\em maximal C$^*-$correspondence of} $X$, for reasons
which will be apparent later.  Clearly $\Li(X)$ may
be rewritten as the closure of:
$$
\left[ \begin{array}{ccl}
\C  & \C X \\
X^* \C & \D  
\end{array}
\right] \; \; .
$$
Example \ref{t2} shows that $C^*_{max}(XX^*) \cong
C^*_{max}(XX^*) \cong C([0,1])$ and $W = 
C_0((0,1])$ if
$X$ is a one dimensional operator space.

Let $H, K$ be general Hilbert spaces.
Note that by the aforementioned modification of
the Lemma 7.1 from \cite{P},
every completely
contractive linear map $T : X \rightarrow B(K,H)$ gives a
completely
contractive unital representation of $\A$ on $H \oplus K$,
and hence a *-representation of 
$\Li(X)$ on $H \oplus K$, and by
restriction, a unital
*-representation $\pi_T$ of $C^*_{max}(X^*X)$
on $K$.   Notice that $\pi_T(x^*y) = T(x)^*T(y)$, for all
$x,y \in X$, and clearly there can be only one such
unital
*-representation of $C^*_{max}(X^*X)$
with this property.
We shall call this construction of $\pi_T$ from $T$
the universal property
of $C^*_{max}(X^*X)$.  One has a
similar universal property
 for $C^*_{max}(XX^*)$.

The converse is also true,
namely, that any unital
*-representation $\pi$ of $C^*_{max}(X^*X)$ on
a Hilbert space $K$
gives rise to a completely contractive
linear map $S_\pi : X \rightarrow B(K,H_\pi)$, for
some Hilbert space $H_\pi$.  Indeed let $M = X \otimes K$,
and define a semi-inner-product on $M$ by
$\langle x_1 \otimes \eta_1 , x_2 \otimes \eta_2 \rangle
= \langle \pi(x_2^* x_1) \eta_1 , \eta_2 \rangle$.
We define $H_\pi$ or $X \otimes_\pi H$
to be the completion of the quotient
of $M$ by the null vectors in $M$.  Then $H_\pi$ is
a Hilbert space, and we define $S_\pi(x)(\eta) = x \otimes
\eta$, for $x \in X, \eta \in K$.  It is easily checked
that $S_\pi$ is completely contractive, and
also that $\pi(x^*y) = S_\pi(x)^*S_\pi(y)$, for all
$x,y \in X$.  Thus $\pi_{S_\pi} = \pi$.  We 
also note that $[S_\pi(X)(K)]$ is dense in $H_\pi$.

Finally note that if one begins with
a completely contractive
linear $T : X \rightarrow B(K,H)$, and then
forms the
associated representation $\pi = \pi_T$
of $C^*_{max}(X^*X)$ on
$K$ as above, and then produces a Hilbert space
$H_\pi$ and complete contraction $S = S_{\pi_T}$ as
in the last paragraph.  Then it is clear that
there is a canonical
(and indeed unique) isometry
$U : H_\pi \rightarrow H$ with the property
that  $U S = T$.  We note that $U$ is a unitary
if and only if 
the span of $T(X)(K)$ is dense in $H$.  The above seems
to be some kind of `polar decomposition' for operator
spaces.  We have written a general $T$ as the
composition of an isometry, and a map $S$ of the 
`standard form' $S(x)(\eta) = x \otimes \eta$.

If one begins with a 
*-representation $\theta$ of C$^*_{max}(XX^*)$ on $H$,
one defines $K_\theta = X^* \otimes_\theta H$ similarly,
and we define a
 map $R_\theta : X \rightarrow B(K_\theta,H)$
given by $R_\theta(x)(y^* \otimes \zeta) = 
\theta(xy^*)(\zeta)$.  Now $R_\theta(x) R_\theta(y)^*
= \theta(xy^*)$.   If $\theta$ comes from a 
map $T : X \rightarrow B(K,H)$ , via the universal
property of C$^*_{max}(XX^*)$,   then one can 
easily check that there is a coisometry $V : K
\rightarrow K_\theta$ such that $T = R_\theta V$.
Also, $V$ is unitary if and only if 
$\cap_{x \in X} \ker T(x) = (0)$.

 The universal property
of $C^*_{max}(X^*X)$ is reminiscent of the property of the
universal C$^*-$algebra C$^*(X)$ of an operator space $X$
\cite{Pes}.  However,
C$^*(X) \cong C^*_{max}(OA(X))$, where
$OA(X)$ is the universal operator algebra of an operator
space discussed in
\cite{Pis}.  Indeed there is no obvious inclusion of
$X$ in $C^*_{max}(X^*X)$.   It is also not true
that $\Li(X)$ coincides with the universal
C$^*-$algebra generated by the operator system
${\mathcal S}$.  To see this, observe that the
latter C$^*-$algebra is shown in \cite{KW} to be
nonexact if ${\mathcal S} = M_2$, whereas in this case
$X = {\mathbb C}$ and $\Li(X)$ is the
(nuclear) C$^*-$algebra in Example \ref{t2}.  However
it is clear that $\Li(X)$ is always
a quotient C$^*-$algebra 
of the C$^*-$algebra of the operator system 
${\mathcal S}$.

It would be interesting to study these universal
C$^*-$algebras
for some of the common finite dimensional 
operator spaces $X$.  Understanding 
$C^*_{max}(X^*X)$ for $X = \ell_n^1$, for example, 
corresponds 
to understanding a certain von Neumann
type inequality.  One must find $n$ universal contractions
which together with their adjoints,
satisfy certain polynomial inequalities.  
As far as we know, this particular
type of von Neumann
type inequality, or such universal C$^*-$algebras,
have not been studied. 

As noted in \cite{Zh}, the subalgebra
$C^*_{0}(XX^*)$ of $C^*_{max}(XX^*)$ generated by $XX^*$ (but
not $e_1$), is,
by construction, strongly Morita equivalent to
$C^*_{0}(X^*X)$.
Thus it is not strange that
these C$^*-$algebras have the same `representation theory'.
We can rephrase the construction of $S_\pi$ from
$\pi$ given above as follows.  A unital
*-representation $\pi$
of $C^*_{max}(X^*X)$ on a Hilbert space $K$,
restricts to
a *-representation of $C^*_{0}(X^*X)$.
By the basic theory of strong Morita equivalence,
$\pi$
gives rise to a canonical second Hilbert space $H$
(which may be obtained as the `interior' or `module Haagerup'
 tensor product of $W$ and $K$ (see \cite{La,Bna} for
example)).  We also obtain a canonical
*-representation of $C^*_{max}(X^*X)$ on $H$, and a canonical
*-representation of $\Li(X)$
on the Hilbert space $H \oplus K$.  This is
the whole point of `induced representations' \cite{Ri2}.
By restriction to the 1-2 corner,
we obtain the canonical completely
contractive linear map $S_\pi : X \rightarrow B(K,H)$.

Looking at $C^*_{max}(X^*X)$ from the point of
view of C$^*-$modules is perhaps the best way to
formulate its universal property.  Namely,
if one takes $W = e_1 \Li(X) e_2$, then
as we just saw, $W$ is a right C$^*-$module over
$C^*_{max}(X^*X)$ , and there is an obvious
complete isometry $i : X \rightarrow W$ such that
the identity and the
range of $i(X)^*i(X)$ generates $C^*_{max}(X^*X)$.
Moreover, any completely contractive map $T : X
\rightarrow Z$ into a right C$^*-$module $Z$ over
$\B$, say, give rise to a (necessarily unique)
unital *-homomorphism $\pi : C^*_{max}(X^*X) \rightarrow \B$
such that $\pi(i(y)^*i(x)) = \langle T(y) \vert T(x)
\rangle$ for all
$x,y \in X$.  This may be seen by applying the
previous universal property of $C^*_{max}(X^*X)$
to $T$, the latter viewed as a map into the range
of a concrete faithful *-representation of
the linking C$^*-$algebra of $Z$.
 Conversely any unital
*-homomorphism $\pi$ from  $C^*_{max}(X^*X)$ into a
C$^*$-algebra $\B$, restricts to a
*-homomorphism $\pi'$ on $C^*_{0}(X^*X)$, which induces
a quotient Morita context for the
range of $\pi'$ as in \cite{Ri3} \S 3.  This is
the `pushout' construction discussed
at the end of our last section.
 We obtain a
C$^*-$bimodule $Z$ and a complete contraction
$T : X \rightarrow Z$  with the
$\pi(i(y)^*i(x)) = \langle T(y) \vert T(x) \rangle$
property.

\vspace{4 mm}

We now generalize the above to give the
 linking algebras, and
the `maximal C$^*-$correspondence' of an operator bimodule.  
Let $X$ be an $\A-\B$-operator module, and suppose
that $\C$ and $\D$ are the maximal C$^*-$algebras of 
$\A$ and $\B$ respectively.  Let $\tilde{X}$ be
the $\C-\D$-dilation
of $X$, namely 
$\tilde{X} = \C \otimes_{h\A} X \otimes_{h\B} \D$.
As we mentioned at the end of
 \S 3, this `bi-dilation' has the
obvious universal property.  For simplicity in what follows
we shall assume $\A$ and $\B$ are unital. The more general
case of c.a.i.'s follows from the unital case by the 
usual tricks.   
By the Christensen-Effros-Sinclair representation
theorem for operator bimodules over C$^*-$algebras, 
there is
a concrete representation of $\C, \D$, and
$\tilde{X}$ on 
Hilbert spaces $H'$ and $K'$ say, in such a way that the 
module
actions become concrete multiplication of operators. 
It can also be arranged, by direct summing with the
identity representations of the C$^*-$algebras,  
that the representations of $\C$ and $\D$ on $H'$
and $K'$ respectively in the
Christensen-Effros-Sinclair representation are 
faithful (and hence completely isometric).
We shall refer to this as a {\em faithful 
CES representation}.  With respect to one fixed 
faithful CES representation we
form an operator system ${\mathcal S}$ in 
$B(H' \oplus K')$ consisting of 
matrices
$$
\left[ \begin{array}{ccl}
c & \tilde{x} \\
\tilde{y}^* & d
\end{array}
\right]
$$
where $c \in \C, d \in \D, \tilde{x}, \tilde{y} 
\in \tilde{X}$.
Paulsen's $2 \times 2$ matrix tricks, used in 
the proof of Lemma
7.1 in \cite{P}, or 
theorem 2.4 in \cite{Su}, for example
may be used in an analogous and
 straightforward fashion to yield the following:

\begin{theorem}  
\label{pr} Given a *-representation 
$\theta$ of $\C$ on a Hilbert space
$H$ and a *-representation $\pi$ of  $\D$
on a Hilbert space
 $K$, and given a completely contractive $\C-\D$-module map 
$\Phi : \tilde{X} \rightarrow B(K,H)$, then the map
$\Psi$ from ${\mathcal S}$ into $B(H \oplus K)$ defined by
$$
\left[ \begin{array}{ccl} c & \tilde{x} \\ 
\tilde{y}^* & d \end{array} 
\right] \; \; \mapsto \; \; \left[ 
\begin{array}{ccl} \theta(c) & \Phi(\tilde{x}) 
\\ \Phi(\tilde{y})^* & \pi(d)
\end{array}
\right]
$$
is completely positive.
\end{theorem}

We will omit the proof.
 From this one can immediately deduce that 
the operator system structure on ${\mathcal S}$ is 
independent of the particular faithful CES representation of 
$\tilde{X}$.   Note that this theorem, seems to give a
direct proof of \ref{wit}, using the idea
in \cite{Su} of extending the c.p. map $\Psi$ and
then using 4.2 in \cite{Wi1} to prove that the 1-2 corner
of the extension is still a bimodule map.  We have not 
checked the details, since it is clearly more
trouble than the proof we gave, and moreover still 
requires fussing with the non-nondegeneracy of the
representations involved.

 We define the  {\em upper
triangular operator algebra} 
or {\em upper linking operator algebra}, 
$\U(X)$ to be the set of 
matrices
$$
\left[ \begin{array}{ccl}
a & x \\
0 &b 
\end{array}
\right]
$$
in ${\mathcal S}$, where $a \in \A, b \in \B, x \in X$.  
We define $\U(\tilde{X})$ similarly,
except that the entries come from $\C, \D$ and
$\tilde{X}$.  By looking at the concrete realization
of ${\mathcal S}$ that we began with, it is easily seen
that with respect to
the natural multiplication on $B(H' \oplus K')$,
the spaces $\U(X)$ and $\U(\tilde{X})$
are operator algebras.  Note that the adjoint
in ${\mathcal S}$ of the copy of $X$ in the 1-2 corner,
is $\bar{X}$, the conjugate operator module of $X$
mentioned at the end of \S 1.   We now form 
$\Li(X) = C^*_{max}(\U(X))$, and call this the
linking C$^*-$algebra of $X$.  We let 
$\Li(\tilde{X})$ be $C^*_{max}(\U(\tilde{X}))$.
It is easily seen that 
$\Li(X)$ has a natural decomposition as $2 \times 2$
matrices, and we define 
$\E$ to be the 1-1 corner 
$e_1 \Li(X) e_1$, $W$ to be the 1-2 corner
$e_1 \Li(X) e_2$, and $\F =
e_2 \Li(X) e_2$ .   Here $e_1$ and
$e_2$  are the
copies of the identities of $\A$ and $\B$.

If one begins with a completely contractive
$\A-\B$-bimodule map $T : X \rightarrow B(K,H)$,
where $H$ and $K$ are $\A$- and $\B$-Hilbert modules 
respectively,
then by the universal property of the dilation $\tilde{X},$
 there is a unique completely contractive
$\C-\D$-bimodule  extension 
$\tilde{T} : \tilde{X} \rightarrow B(K,H)$.
By \ref{pr}, we obtain a completely positive unital 
$\Psi : {\mathcal S} \rightarrow B(H \oplus K)$, 
and by restriction, a completely contractive unital
homomorphism $\sigma$ on $\U(\tilde{X})$,  and 
another, $\sigma'$, on $\U(X)$.  Conversely 
any completely contractive unital
homomorphism $\sigma' : \U(X) \rightarrow B(N)$ determines 
a decomposition $N = H \oplus K$, and *-representations
$\pi$ and $\theta$ on $K$ and $H$ respectively, and 
a completely contractive
$\A-\B$-bimodule map $X \rightarrow B(K,H)$.
Thus, as before, there are 1-1
correspondences between the four classes of maps
whose elements we have labeled above with 
symbols $T, \tilde{T} , \sigma , $ and $\sigma'$.

By the universal property of the maximal C$^*-$algebra,
$\sigma$ extends to a *-representation $\tilde{\sigma}$ of
$\Li(\tilde{X})$ .  The restriction of 
$\tilde{\sigma}$ to $\U(X)$ clearly coincides with
$\sigma'$, which shows that $\Li(X)$ may be taken to 
be the C$^*-$algebra inside $\Li(\tilde{X})$ generated
by $\U(X)$.  Hence $\Li(X) = \Li(\tilde{X})$.   
Clearly, one sees also  that ${\mathcal S}$ sits
naturally inside $\Li(X)$.  Thus inside $\Li(X)$
the product $[\C X \D]^{\bar{}}$ is 
completely isometrically isomorphic to
$\C \otimes_{h\A}
X \otimes_{h\B} D$ .
And we can add to our list of 1-1 correspondences
between classes of maps, the correspondence between
completely contractive
$\A-\B$-bimodule maps $T : X \rightarrow B(K,H)$,
and unital *-representations
of $\Li(X)$ on $N (= H \oplus K)$.
  
Notice the above gives a simple way of writing any
completely contractive
$\A-\B$-bimodule map $T : X \rightarrow B(K,H)$
as $P_H \pi(\cdot) \vert_{K}$ where $\pi$
 is a *-representation on $H \oplus K$
of a C$^*-$algebra which contains $X$.  Using the 
CES representation theorem, one can replace the
$B(K,H)$ in the last sentence
by any $\A-\B$-operator bimodule $V$, to represent
an $\A-\B$-bimodule map $T : X \rightarrow V$
as $R \pi(\cdot) V$ where $R, V$ are isometric
or coisometric module maps. 

One may proceed as before, to show that
there are 1-1 correspondences between 
completely contractive
$\A-\B$-bimodule maps $T : X \rightarrow B(K,H)$,
and *-representations $\pi$ and $\theta$ of $\F$ and
$\E$ on $H$ and $K$ respectively.  This all goes through
with no essential changes.  
Again it is interesting that one can write
a general such $T$, as a standard form $S(x)(\eta)
= x \otimes \eta$, multiplied by an isometry
(or coisometry in the other case).  

This universal property of $\F$ (and similarly the 
one for $\E$)
is probably again
best described in terms of C$^*-$modules.
Again $W$ is a right C$^*-$module over $\F$, 
and there is an obvious $\A-\B$-module map $i : 
X \rightarrow W$.
Recall that a right $\C-\G$-C$^*-$correspondence
$V$ (also known as a $\G$-rigged $\C$-module)
is a right 
C$^*-$module over a C$^*-$algebra $\G$, which is
 also a nondegenerate
left Banach $\C-$module (see \cite{Bna} \S 4).
We define a {\em rigging map} to be
a completely contractive
left $\A$-module map $T : X \rightarrow V$ into
a right $\C-\G$-correspondence,   
for which there exists a 
unital completely contractive
homomorphism $\pi : \B \rightarrow
\G$ such that $T(x b) = T(x) \pi(b)$ for all 
$b \in \B$.  Notice that this makes $T$ an
$\A-\B$-bimodule map.  Also
note that the $W$ above is a 
$\C-\D-$C$^*-$correspondence, and that 
the $i : X \rightarrow W$ is a rigging map.
As before, we have a pushout construction.
The
relation 
$\langle T(x) \vert T(y) \rangle =
\pi(i(x)^* i(y))$ for all
$x , y \in X$, determines a correspondence
between 1) rigging maps 
$T : X \rightarrow V$ into a right 
$\C-\G$-correspondence $V$, 
and 2) unital 
*-homomorphisms $\pi'$ from $\F$ into a C$^*$-algebra.
We omit the details, and the standard adaption to the
nonunital case.  The other universal properties of
$\U(X), \Li(X)$ may also be stated in terms of 
C$^*-$modules and C$^*-$correspondences, but we will not
take the time to do that here.

\vspace{4 mm}

\noindent {\bf Final remark.}  We feel that there is
some aspect missing in our   
understanding of operator modules. 
The fact that the notion we called
C$^*-$restrictability  in \cite{BlFTM2} is automatic 
suggests strongly the need for a good test for an
$\A$-operator module to be a $\C$-operator module.

We thank Christian Le Merdy for many 
useful comments on previous versions of this paper,
and for answering several questions.

\setcounter{section}{3}

\end{document}